# A More Robust Multi-Parameter Conformal Mapping Method for Geometry Generation of any Arbitrary Ship Section

Mohammad Salehi [a], Parviz Ghadimi [b*], Ali Bakhshandeh Rostami [c]

[a] Graduate Student at Dept. of Marine technology, Amirkabir University of Technology, Teheran, Iran

[b] Associate Professor at Dept. of Marine technology, Amirkabir University of Technology, Teheran, Iran

[c] M.Sc. Graduate at Dept. of Marine technology, Amirkabir University of Technology, Teheran, Iran

## Abstract

The central problem of strip theory is the calculation of potential flow around 2D sections. One particular method of solutions for this problem is conformal mapping of the body section to the unit circle over which solution of potential flow is available. Here, a new multi-parameter conformal mapping method is presented that can map any arbitrary section onto a unit circle with good accuracy. The procedure for finding the corresponding mapping coefficients is iterative. The suggested mapping technique is proved able to map any chine, bulbous, large and fine sections, appropriately. Several examples of mapping the symmetrical as well as the nonsymmetrical sections have been demonstrated. For the symmetrical and nonsymmetrical sections, the results of the current method are compared against other mapping techniques and the currently produced geometries display good agreement with the actual geometries.

## 1. Introduction

Dynamic analysis of floating and sailing bodies is an important part of basic design of marine structures that helps the designers to estimate the dimensions and forms of the bodies. Rigid body motions are defined by six equations and each equation belongs to one degree of freedom. To obtain all body motions, these six equations are coupled and solved simultaneously. These differential equations consist of mass, damping, stiffness terms and wave loads. Added mass or added inertia are taken into account and is due to the acceleration of fluid that surrounds the body, while the damping term is due to the reduction in the energy level of the body which generates the free surface waves during the body motions. In order to compute the added mass and added inertia, damping and wave loads, it is necessary to solve the potential flow problem around the body which must satisfy some boundary conditions such as the body condition, free surface condition and radiation condition.

One method of solution for calculating the three dimensional problems is to break down the geometry into some two dimensional sections, solve each section separately, and subsequently superpose the sectional results to gain the response of the main geometry. This method is called strip theory whose simplicity and high accuracy has caused its extensive application in calculation of the potential flows. Descriptions of the background of strip theory are provided by Salvesen et al [1], who also developed a detailed formulation of the theory. Their version has been the most widely used strip theories. One particular scheme for application of the strip theory in solving the potential flow problems is conformal mapping. This is accomplished by mapping any arbitrary shape to a simpler form like a circle for which there is an analytical solution available for potential flow. By using this solution and the inverse of the applied conformal mapping, the potential flow around the intended arbitrary shape is determined. Ursell [2] derived the potential flow around a circular cylinder which involved the heaving in the free surface. This potential flow satisfied the free-surface condition as well as the body boundary condition. Therefore, by using the Ursel method to calculate the potential flow around a cylinder with an arbitrary section, one has to find the relation between a circle and the arbitrary section. The relation between the arbitrary section and the circle is established by some number of coefficients and parameters.

There are other types of conformal mapping that are used for flow analyses in the field of hydrodynamics. Schwartz-Christoffel method is one mapping scheme that has the superior benefits of mapping the sections with high curvature, like the chine sections, the rise of floors and the wedge shapes. This method maps the section to a line or a polygon. Ghadimi *et al* [3] used Schwartz- Christoffel method to map the wedge section onto a line for solving the water entry problem. They calculated the pressure distribution for several entry angles in potential flow. Hasheminejad and Mohammadi [4] mapped the elliptical and circular domains to rectangular domain to study the effect of anti-slosh baffles on free liquid oscillations in partially filled horizontal circular and elliptical tanks. Modified Schwarz-Christoffel mappings, using approximate curve factors, were also presented by Andersson [5] that was able to map the upper-half plane to a polygonal region.

A mapping technique that involves the transformation of any arbitrary section to a unit circle is of immense importance. One good example of this type of mapping scheme is the Riemann mapping theory that was illustrated by Papamichael and Stylianopoulos [6]. In this method, three real conditions must be imposed in order to make the conformal mapping unique. In other words, the problem of determining the mapping in this method has three degrees of freedom. Another example of transforming a rectangle to a unit circle is called the KT-TG method. For a rectangular section, the conformal transformations combine successively the Karmann–Trefftz (KT) transformation (see Halsey [7], for example) which removes the corners and the Theodorsen–Garrick (TG) [8] transformation which turns the intermediate near circle into a perfect circle. This method was used for computing the hydrodynamic forces acting on bodies of arbitrary shape in viscous flow by Scolan and Etienne [9]. Carmona and Fedorovskiy [10] studied some properties of Carathéodory domains and their conformal mappings onto the unit disk. Brown [11] presented a two parameter conformal mapping from a disk to a circular-arc quadrilateral, symmetric with respect to the coordinate axes. Brown and Porter [12] developed a conformal mapping from the circular-arc quadrilaterals with four right angles onto a unit disk. Lui et al [13] proposed a representation of general 2D domains with arbitrary topologies using conformal geometry. A natural metric can be defined on the proposed representation space, which gives a metric to measure dissimilarities between objects. The main idea was to map the exterior and interior of the domain conformally onto unit disks and circular domains. As a general comment, although the use of Green functions with or without forward speed, produces more accurate results for the ship hull forms, but due to their simplicity, conformal mapping techniques are still suitable for seakeeping calculations.

As mentioned earlier, there are many types of conformal mappings that are used for the analysis of flow around marine structures, but the most important and interesting mapping technique is that which can map any ship sections onto a unit circle, which is the focus of the current study. In this regard, the first and most famous conformal mapping method is two-parameter Lewis conformal mapping which was presented by Lewis [14]. The base of this method was the similarity of the breadth and draft and area of the section and those of the mapped section. This method is able to map a circle to a conventionally arbitrary shape but does not produce good results for any arbitrary sections. Thereafter, three-parameter conformal mapping was introduced by Tasai [15]. They extended two-parameter Lewis conformal mapping to three-parameter conformal mapping, by taking into account the first order moments of half of the real body cross sectional surface and that of the mapped section about the x and y-axes.

In 1969, von Kerczek and Tuck [16] for the first time presented multi-parameter conformal mapping based on Least Square Method by using an iterative procedure to solve nonlinear system of equation which resulted in finding the mapping coefficients. This method was later simplified by equating the vertical position of the centroid of the real body cross section and that of the mapped section by Athanassoulis and Loukakis [17]. The result of this method was proved to be better than that of Lewis method, but this method was still not able to map any or all sections. In the meantime, Anghel and Ciobanu [18] showed how to extend the Lewis transformation for obtaining the contour of the ship's cross section of different types of ships. Moreover, they demonstrated how a more accurate transformation of the cross sectional hull shape can be obtained by using a larger number of parameters. Subsequently, a multi-parameter conformal mapping was presented by Tasai [19, 20] which was further developed by de Jong [21].

Rosén and Palmquist [22] used the de Jong and Lewis mapping method to compute the added mass and damping coefficients of a ship, for studying the parametric rolling in waves. Westlake and Wilson (WW) [23] developed a technique for carrying the multi-parameter conformal mapping of the ship sections onto a unit circle. The basic idea was to calculate one unknown mapping parameter after another when the starting known mapping coefficient was calculated for N=0, where N was the number of mapping coefficients. Also, in their mapping method, the corresponding angle in the unit circle is calculated by interpolation, using the line length method corresponding to each defined point. Another multi-parameter conformal mapping is the method that was presented by Journée and Adegeest [24], but they only studied the symmetrical sections.

In the currently proposed multi-parameter mapping method, the least square scheme that was presented by de Jong [21], has been used to compute the mapping coefficients, but a new technique has been introduced that is drastically different in detail and can map any arbitrary symmetrical and non–symmetrical sections onto a unit circle in a much quicker way. However, contrary to the method suggested by Westlake and Wilson [23], the value of the first mapping coefficients based on de Jong's [21] general mapping equations, is assumed to be one, all the mapping parameters are calculated simultaneously, and by making a simple assumption, an equation is obtained and solved for the corresponding angles in the unit circle plane. Furthermore, in the proposed method, the initial guess for the mapping coefficients for the non-symmetrical sections, will be assumed as the average of the two mapping coefficients found based on the approach adopted for the symmetrical sections and for the right and left (half) distance from the centerline on the waterline. The newly introduced schemes make the proposed multi-parameter mapping technique more robust, shorten the process of mapping and decrease the computing time.

To recapitulate the main ideas of this paper, a more robust multi-parameter conformal mapping method is presented that is able to map any arbitrary sections onto a unit circle. This method, shortens the process of mapping, decreases the computing time, and has the ability to map the extreme bulbous sections, the non-symmetrical sections, the large area coefficient sections, sections with hard chine and fine sections among others. The process of mapping is outlined and detailed clearly and many practical examples of the proposed method have been offered which demonstrate its great accuracy. Results of the current mapping technique, in the case of symmetrical sections, have been plotted and compared against the results of the Lewis mapping method which indicate better agreement of the current results with the points on the real sections. To further indicate the abilities of the current mapping method in the cases of rectangular, large bulbous and fine sections for both symmetrical and non-symmetrical form, comparisons are also made between the current mapping results and the results of Westlake and Wilson [23] Mapping scheme. These comparisons demonstrate that in symmetrical cases, the current mapping method is more robust while in non-symmetrical cases, the results are very similar to that by Westlake and Wilson [23] mapping technique. For a more clear observation of the relations between the real section points and their corresponding angles in the mapped unit circle plane, corresponding angles for the case of symmetric rectangular section has been displayed in a figure.

## 2. Basic formulation and assumptions for conformal mapping

The tranformation of any arbitrary section (in Cartesian plane) to a unit circle plane, using a multi-parameter conformal mapping is considered to be [21]:

$$z = F \sum_{n=0}^{N} a_{2n-1} \xi^{-(2n-1)} \tag{1}$$

where $z = x + iy$ represents a Cartesian plane (as shown in Fig.1), $x$ is the real part of $z$, $y$ is the imaginary part of $z$, $\xi = ie^{\beta} e^{-i\theta}$ represents a unit circle, $F$ is the scale factor, $a_{-1} = +1$ and $a_{2n-1}(n=1,...,N)$ are the mapping coefficients, $\theta$ is the corresponding angle in the mapped unit circle plane, $\beta$ is considered another scaling factor indicating how close any mapped section is from the mapped ship section, and $N$ is the number of mapping parameters.

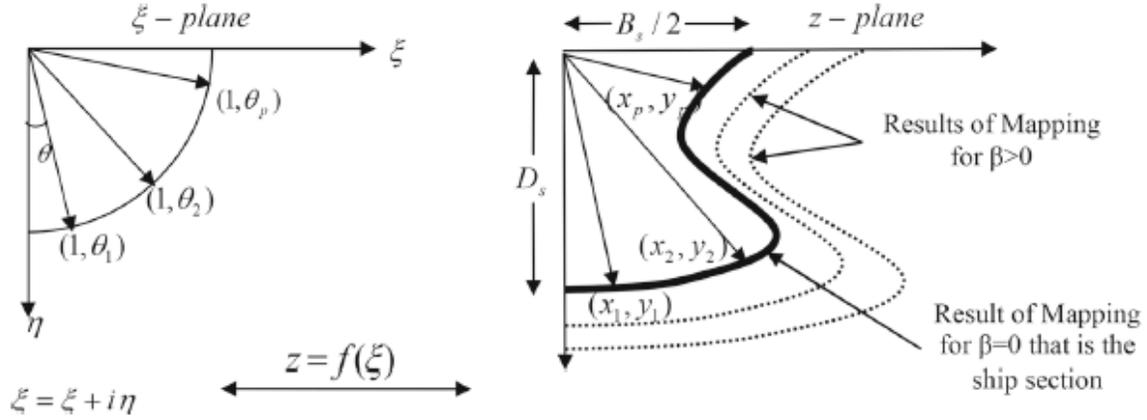

Fig.1. Actual coordinate system for the arbitrary section (right) and the transformed unit circle plane (left).

Based on the diagrams presented in Fig.1, the $x$ and $y$ coordinates of the ship section in Eq.1 can be rewritten in the form of

$$x = -F \sum_{n=0}^{N} (-1)^n a_{2n-1} e^{-(2n-1)\beta} \sin((2n-1)\theta) \tag{2}$$

$$y = +F \sum_{n=0}^{N} (-1)^n a_{2n-1} e^{-(2n-1)\beta} \cos((2n-1)\theta) \tag{3}$$

As shown in Fig.1, $\beta=0$ results in the mapping of the boundary of ship section in which case, Eqs.2 and Eq.3 become

$$x_0 = -F \sum_{n=0}^{N} (-1)^n a_{2n-1} \sin((2n-1)\theta) \tag{4}$$

$$y_0 = +F\sum_{n=0}^{N}(-1)^n a_{2n-1} \cos((2n-1)\theta) \tag{5}$$

Draft and breadth of each section are defined as follows:

$$D_s = F\sigma_b \quad \text{with} \quad \sigma_b = \sum_{n=0}^{N}(-1)^n a_{2n-1} \tag{6}$$

$$B_s = 2F\sigma_a \quad \text{with} \quad \sigma_a = \sum_{n=0}^{N} a_{2n-1} \quad , \quad F = \frac{B_s}{2\sigma_a} \tag{7}$$

Usually, the arbitrary sections are defined by some discrete points such as $(x_i, y_i)(i = 0,1,2,...,I)$. There is also an angle $\theta_i$ in each unit circle plane corresponding to each point of the arbitrary section. By applying the appropriate mapping coefficients and the corresponding values of $\theta_i$ in Eqs.4 and 5, points on the arbitrary section will be produced. Therefore, the conformal mapping problem reduces to that of finding the mapping coefficient and the values of $\theta_i$. The number of mapping coefficients, i.e. N, is an arbitrary value.

To determine the exact breadth at the water line and the draught, there are two points for which their corresponding angles in the unit circle are known. These two points are found from Eqs.8 and 9 and are given below.

$$x = \frac{B_s}{2} \quad , \quad y = 0 \quad\quad \theta = \frac{\pi}{2} \tag{8}$$

$$x = 0 \quad , \quad y = D_s \quad\quad \theta = 0 \tag{9}$$

The method starts by guessing the initial values for mapping coefficients and subsequently using these coefficients to find the corresponding angles $(\theta_i)$ in the unit circle for each point of the arbitrary section. Next, the new values of the mapping coefficients are found by using the computed angles $(\theta_i)$. This computation is continued until the summation of squared errors between the mapped and exact points of the section becomes less than a designated error tolerance value.

## 3. The Proposed Multi-Parameter Conformal Mapping Technique

This method can be described and utilized for both symmetrical and nonsymmetrical sections. In this section, initially the symmetrical section conformal mapping and subsequently the conformal mapping for non-symmetrical sections will be described.

### 3-1. Conformal Mapping of the Symmetrical sections

The results of conformal mapping of $(x_i, y_i)$ are assumed to be $(x_{0i}, y_{0i})(i = 0,1,2,...,I)$. First, by using the two-parameter Lewis conformal mapping method, the initial values of conformal mapping coefficients are defined. Then, for minimization of the square of differences between $(x_i, y_i)$ and $(x_{0i}, y_{0i})(i = 0,1,2,...,I)$, two conditions or constraints should be satisfied which are as follows:

a) **First condition or constraint:**

Points $(x_{0i}, y_{0i})$ of the transformed contour lie on the normal of the arbitrary section. This normal passes through the points $(x_i, y_i)$ of the section contour, as shown in Fig.2. By satisfying this condition, the unknown values of $\theta_i$ $(i = 0,1,2,...,I)$ will be obtained.

b) **Second condition or constraint:**

The summation of square of differences of $(x_i, y_i)$ and $(x_{0i}, y_{0i})(i = 0,1,2,...,I)$ must be less than an arbitrary designated value. By satisfying this condition and utilizing the obtained values of $\theta_i$ found from the first condition, the values of mapping coefficients will be produced. These determined coefficients are used as initial values of mapping coefficients. Computing the unknown values of $\theta_i$ is repeated with these new initial mapping coefficients and the mapping coefficients are updated by using the obtained values of $\theta_i$. This procedure is repeated until the summation of square of differences of $(x_i, y_i)$ and $(x_{0i}, y_{0i})$ becomes less than an arbitrary designated error tolerance value.

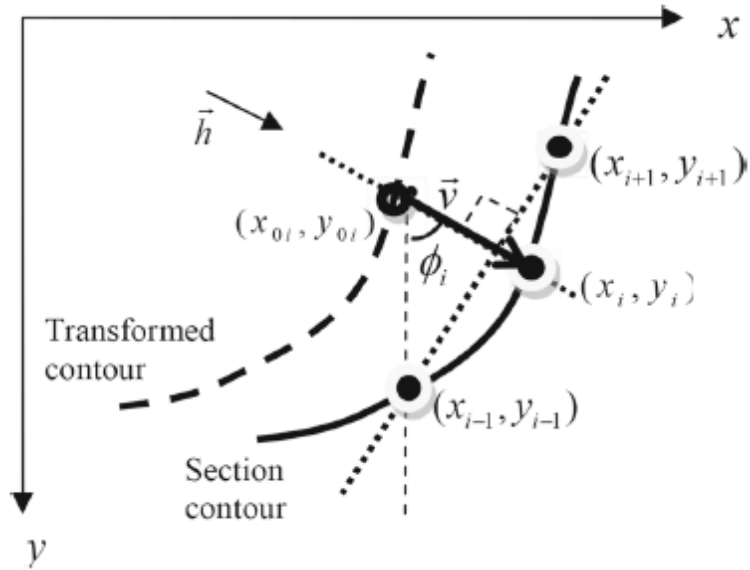

Fig.2. Normal line to the arbitrary section that passes through $(x_i, y_i)$ and $(x_{0i}, y_{0i})$.

In order to satisfy the first condition, the vector passing through the mapped point of $(x_{0i}, y_{0i})$ and normal to the secant line connecting two points of $(x_{i-1}, y_{i-1})$ and $(x_{i+1}, y_{i+1})$, i.e. $\vec{h}$, must be parallel to the vector passing through the points $(x_{0i}, y_{0i})$ and $(x_i, y_i)$, i.e. $\vec{v}$. Therefore, their cross product must be zero.

The unit normal vector passing through the mapped point of $(x_{0i}, y_{0i})$ and normal to the secant line connecting two points of $(x_{i-1}, y_{i-1})$ and $(x_{i+1}, y_{i+1})$ on the real section is defined as

$$\vec{h} = \sin\phi_i i + \cos\phi_i j \tag{10}$$

where $\phi_i$ is the angle between the vertical line that passes through point $(x_{0i}, y_{0i})$ and the unit normal vector. The angle $\phi_i$ is constant throughout the iteration because the normal line is constant, i.e. $\vec{h}$ is constant, but angle $\theta_i$ varies in each iteration and converges to the exact corresponding angle in the unit circle plane. The cause of this variation is the change in vector $\vec{v}$.

The vector connecting the two points of $(x_i, y_i)$ and $(x_{0i}, y_{0i})$ is given as

$$\vec{v} = (x_i - x_{0i})i + (y_i - y_{0i})j \tag{11}$$

As mentioned earlier, the two vectors in Eqs. (10) and (11) are parallel. Accordingly, their cross product is zero which leads to the following equation:

$$(x_i - x_{0i})\cos\phi_i - (y_i - y_{0i})\sin\phi_i = 0 \tag{12}$$

Based on Fig.2, $\phi_i$ can be found from the equations (13) and (14) given as:

$$\cos\phi_i = \frac{x_{i+1} - x_{i-1}}{\sqrt{(x_{i+1} - x_{i-1})^2 + (y_{i+1} - y_{i-1})^2}} \tag{13}$$

$$\sin\phi_i = \frac{-y_{i+1} + y_{i-1}}{\sqrt{(x_{i+1} - x_{i-1})^2 + (y_{i+1} - y_{i-1})^2}} \tag{14}$$

Substituting Eqs.4, 5, 13, 14 into Eq.12, would yield in

$$x_i \cos\phi_i + F \cos\phi_i (\sum_{n=0}^{N} (-1)^n a_{2n-1} \sin(2n-1)\theta) - y_i \sin\phi_i + F \sin\phi_i (\sum_{n=0}^{N} (-1)^n a_{2n-1} \cos(2n-1)\theta) = 0 \tag{15}$$

Solving Eq.(15) for each point of the arbitrary section, i.e. $(x_i, y_i)(i = 0,1,2,...,I)$, results in the corresponding values of $\theta_i$ in the unit circle plane. This nonlinear equation is solved by implementing the bisection method.

Applying the second condition, would result in

$$E = \sum_{i=0}^{I} e_i < \sigma_E \tag{16}$$

Where $\sigma_E$ Is designated error tolerance value

And $e_i$ is

$$e_i = (x_i - x_{0i})^2 + (y_i - y_{0i})^2 \tag{17}$$

Substituting Eqs.4 and 5 in Eq.17, would yield in

$$E = \sum_{i=0}^{I} \left\{ \left( x_i + \sum_{n=0}^{N} (-1)^n Fa_{2n-1} \sin((2n-1)\theta_i) \right)^2 + \left( y_i - \sum_{n=0}^{N} (-1)^n Fa_{2n-1} \cos((2n-1)\theta_i) \right)^2 \right\} \tag{18}$$

Therefore, the new values of $Fa_{2n-1}$ should be computed in such a way that the value of $E$, by applying least square method, becomes minimum as in

$$\frac{\partial E}{\partial (Fa_{2j-1})} = 0 \text{ for } j = 0,...,N. \tag{19}$$

This minimization results in $N+1$ equations as in

$$\sum_{n=0}^{N}\sum_{i=0}^{I} (-1)^n Fa_{2n-1} \cos((2j-2n)\theta_i) = \\ \sum_{i=0}^{I} (-x_i \sin((2j-1)\theta_i) + y_i \cos((2j-1)\theta_i)) \quad \text{For: } j=0,...N. \tag{20}$$

For a better and faster convergence in solving the system and obtaining the exact breadth and draft, the last two equations in the system of equations (i.e. for $j=N-1$ and $j=N$ in Eq.20) are replaced by the related equations of breadth and draft as in Eqs.22 and 23. Consequently, Eq.20 can now be written as Eqs.21, 22, and 23 as follows:

$$\sum_{n=0}^{N}\sum_{i=0}^{I} (-1)^n Fa_{2n-1} \cos((2j-2n)\theta_i) = \\ \sum_{i=0}^{I} (-x_i \sin((2j-1)\theta_i) + y_i \cos((2j-1)\theta_i)) \quad \text{For: } j=0,...N-2 \tag{21}$$

$$\sum_{n=0}^{N}(-1)^n Fa_{2n-1} = D_s \tag{22}$$

$$\sum_{n=0}^{N} Fa_{2n-1} = \frac{B_s}{2} \tag{23}$$

This system of equations can be rewritten in matrix form AX=B as in

$$A = \begin{pmatrix} I+1 & -\sum_{i=0}^{I}\cos 2\theta_i & -\sum_{i=0}^{I}\cos 4\theta_i & -\sum_{i=0}^{I}\cos 6\theta_i & \cdots & (-1)^N \sum_{i=0}^{I}\cos 2N\theta_i \\ \sum_{i=0}^{I}\cos 2\theta_i & -(I+1) & \sum_{i=0}^{I}\cos 2\theta_i & -\sum_{i=0}^{I}\cos 4\theta_i & \cdots & (-1)^N \sum_{i=0}^{I}\cos(2-2N)\theta_i \\ \sum_{i=0}^{I}\cos 4\theta_i & -\sum_{i=0}^{I}\cos 2\theta_i & (I+1) & -\sum_{i=0}^{I}\cos 2\theta_i & \cdots & (-1)^N \sum_{i=0}^{I}\cos(4-2N)\theta_i \\ \sum_{i=0}^{I}\cos 6\theta_i & -\sum_{i=0}^{I}\cos 4\theta_i & \sum_{i=0}^{I}\cos 2\theta_i & -(I+1) & \cdots & (-1)^N \sum_{i=0}^{I}\cos(6-2N)\theta_i \\ \vdots & \vdots & \vdots & \vdots & \vdots & \vdots \\ \sum_{i=0}^{I}\cos(2(N-2))\theta_i & -\sum_{i=0}^{I}\cos(2(N-2)-2)\theta_i & \sum_{i=0}^{I}\cos(2(N-2)-4)\theta_i & \sum_{i=0}^{I}\cos(2(N-2)-6)\theta_i & \cdots & (-1)^N \sum_{i=0}^{I}\cos(2(N-2)-2N)\theta_i \\ 1 & -1 & 1 & -1 & \cdots & (-1)^N \\ 1 & 1 & 1 & 1 & \cdots & 1 \end{pmatrix} \tag{24}$$

$$X = \begin{pmatrix} Fa_{-1} \\ Fa_1 \\ Fa_3 \\ Fa_5 \\ \vdots \\ Fa_{2N-1} \end{pmatrix}, \tag{25}$$

and

$$B = \begin{pmatrix} \sum_{i=0}^{I} \{x_i \sin\theta_i + y_i \cos\theta_i\} \\ \sum_{i=0}^{I} \{-x_i \sin\theta_i + y_i \cos\theta_i\} \\ \sum_{i=0}^{I} \{-x_i \sin 3\theta_i + y_i \cos 3\theta_i\} \\ \sum_{i=0}^{I} \{-x_i \sin 5\theta_i + y_i \cos 5\theta_i\} \\ \vdots \\ D_s \\ \dfrac{B_s}{2} \end{pmatrix} \quad (26)$$

The new values of $(Fa_{2n-1}; n = 0,1,...,N)$ will be computed by solving the system of equations in (21), (22) and (23). The system of equations can be solved by using any direct or iterative numerical methods. In this article, the solution has been performed using the LU decomposition method.

The new values of $\theta_i$ can be obtained by using the new values of $(Fa_{2n-1}; n = 0,1,...N)$ in Eq.15. Applying the newly obtained values of $\theta_i (i = 0,1,2,...,I)$ in Eqs.21, 22, and 23 and solving the system of equations provides new values of $(Fa_{2n-1}; n = 0,1,...N)$. This repeating procedure is continued until the value of $E$ that is computed from Eq.18 becomes less than a designated value.

Finally, since $a_{-1} = +1$, the value of the scale factor $F$ equals to $Fa_{-1}$. Other mapping coefficients will be obtained by dividing the computed values of $(Fa_{2n-1}; n = 0,1,...N)$ by $F$.

### 3-2. Conformal Mapping of the Non-Symmetrical sections

The procedure for conformal mapping of the non-symmetrical sections is similar to the symmetrical sections mapping, with the exceptions of the following three tasks:

1. The corresponding angles in the unit circle plane for two points of the section that lie on the waterline must be found, whereas in the case of symmetrical sections they were known to be 0 and $\frac{\pi}{2}$.

Based on Fig.3, the angle $\phi_i$ is defined by its sine and cosine as

$$\cos\phi_i = \frac{x_i - x_{i-1}}{\sqrt{(x_i - x_{i-1})^2 + (y_i - y_{i-1})^2}}, \qquad (27)$$

$$\sin\phi_i = \frac{-y_i + y_{i-1}}{\sqrt{(x_i - x_{i-1})^2 + (y_i - y_{i-1})^2}}, \qquad (28)$$

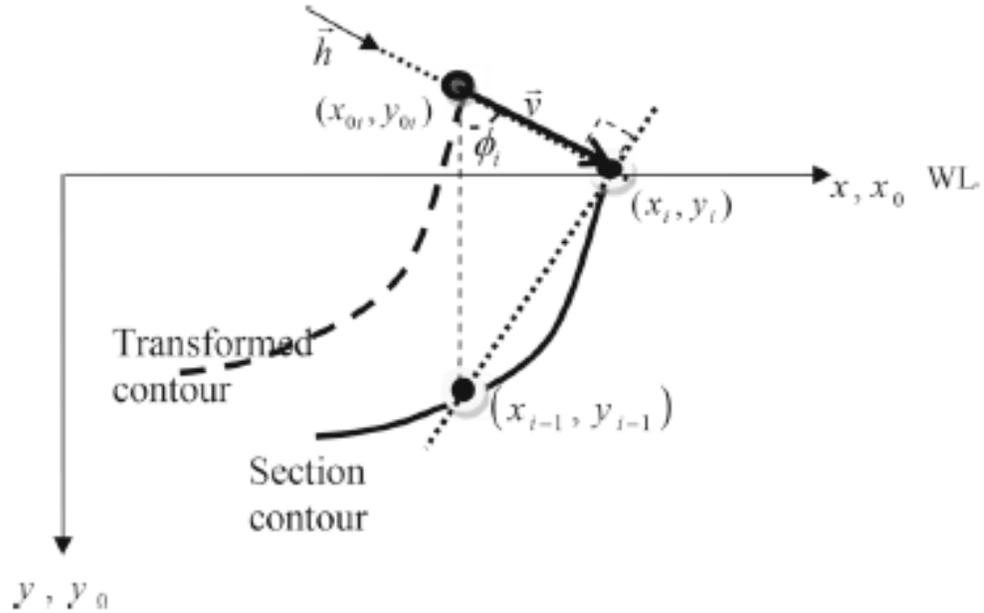

Fig.3. Normal line to the section at the intersection with the waterline(WL).

Therefore, the two values of $\theta$, mentioned above, can be obtained by placing the results of Eqs.27 and 28 into Eq.15 and solving the resulting equation.

2. As mentioned before, for determining the exact breadth and draft in the mapping of the symmetrical sections, two conditions are imposed in the form of two equations on the system of equations (21), but in the non-symmetrical section, this has been proved impractical because solution convergence is not achieved, as solution always diverges. This omission of equations may lead to some error, but this discrepancy only exists between the two points on the waterline and their corresponding points on the mapped section. However, it has been proved, based on the results of the examples provided in this paper, that the reported error is so minimal that it is definitely ignorable.

Therefore, for computing the new values of ($Fa_{2n-1}$; $n=0,1,...N$), the previously computed values of $\theta_i$ should be used in the system of equations (26) as in

$$\sum_{n=0}^{N}\sum_{i=0}^{I}(-1)^n Fa_{2n-1} \cos((2j-2n)\theta_i) =$$
$$\sum_{i=0}^{I}(-x_i \sin((2j-1)\theta_i) + y_i \cos((2j-1)\theta_i)) \qquad \text{For: } j=0,...N \qquad (29)$$

3. Also, as opposed to the conformal mapping of the symmetrical sections for which an initial guess for the mapping coefficients was produced by the Lewis method, for the case of non-symmetrical mapping, the initial guess is obtained in a more complex procedure:

   a. First, by following the same procedure for the symmetrical sections and assuming that its left (half) distance from the centerline on the waterline is ($B_{sl}$), we find the mapping coefficients.
   b. Then, by following the same procedure for the symmetrical sections and assuming that its right (half) distance from the centerline on the waterline is ($B_{sr}$), we find the mapping coefficients.
   c. Finally, the initial guess for the mapping coefficients of the non-symmetrical sections will be assumed as the average of the two mapping coefficients found in parts (a) and (b).

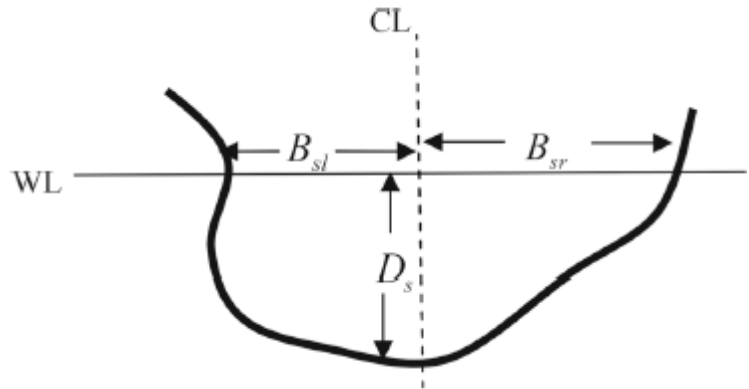

Fig.4. Schematic of a non-symmetrical section.

In the process of mapping the nonsymmetrical sections at the stage of computing the $\theta_i$ values, it is possible that Eq.15 doesn't have a solution. Therefore, in order to avoid this difficulty and to continue with the computations, the value of $\theta_i$ can be assumed in the form of Eq.30.

$$\theta_i = \theta_{i-1} + (\theta_{i-1} - \theta_{i-2}) \tag{30}$$

For better understanding of the current method, a flowchart of the mapping procedure has been plotted in Fig.5.

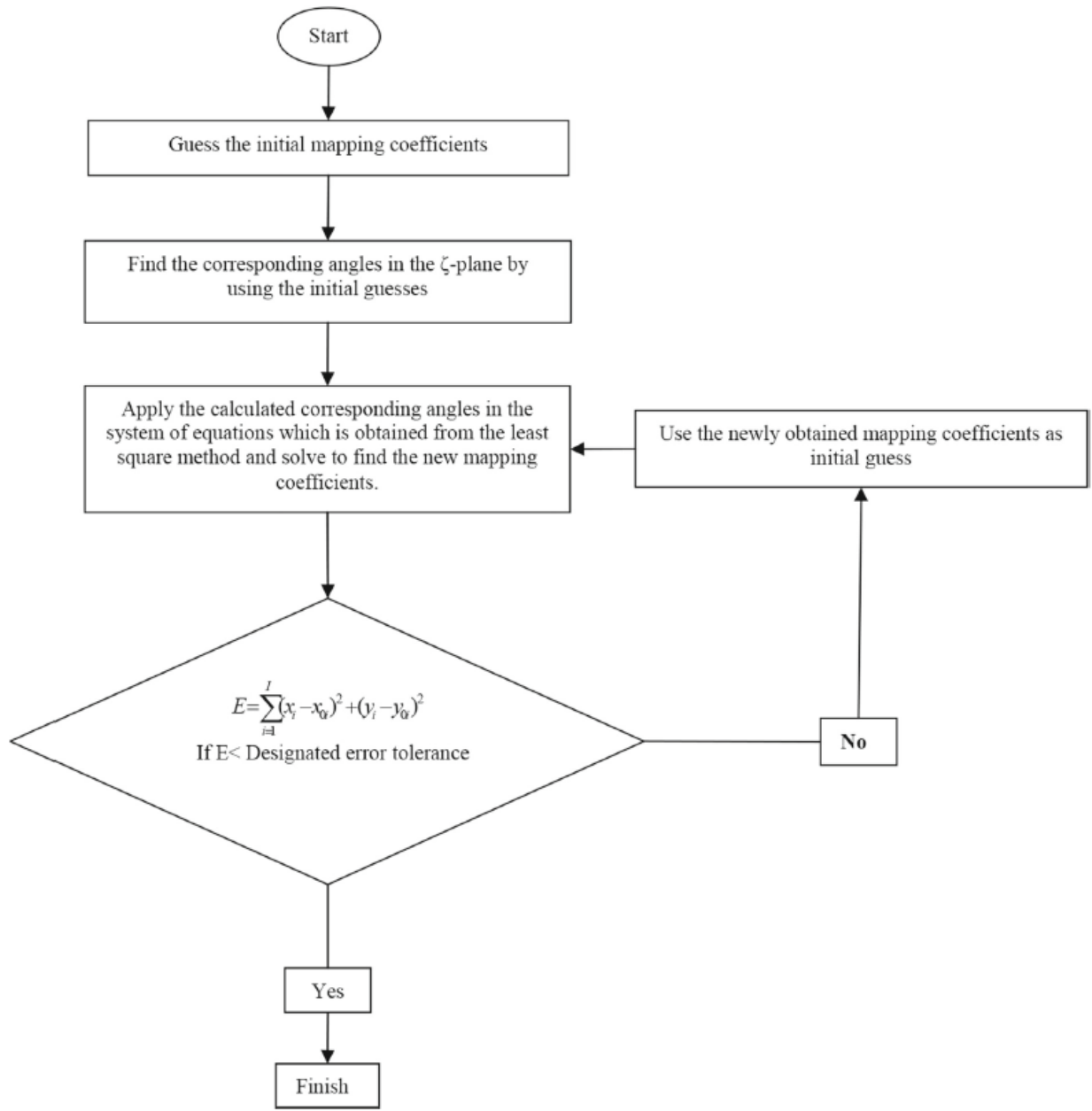

Fig.5 The mapping procedure flowchart.

## 3-3. Finding the Optimum Values of Error (E) and Number of Coefficients (N)

In the mapping process, based on a designated value of N such as N=5, either the program (shown in Fig.5) after a few limited iterations, reaches an optimum E (designated a-priori) or it continues for a large number of iterations without reaching a desirable E. Therefore, a maximum number of 200 iterations is allowed in the mapping process for the symmetric sections and a maximum number of 300 iterations is assigned for the non-symmetric sections in the current computations. Hence, for every value of N, there is a particular value of $E_{min}$ for which the program does not have the ability of reaching a lesser error, anymore. The minimum value of $E_{min}$ found for N=5 is used for setting the designated error tolerance value for N=6 and the process continues. To summarize the complete process of finding the optimum value of E leading to an optimum value of N, the following steps are pursued:

1) We start with N=5 and E=10 as designated error tolerance value.
2) Execute the program according to the flowchart given in Fig.5.
    a. If it does not reach the designated error tolerance of E=10, we move on to N=6 and E=10 as designated error tolerance value.
    b. If it reaches the designated error tolerance of E=10 before 200 iterations for symmetrical cases and 300 iterations for nonsymmetrical cases; then, a program which is written for finding an optimum error for any N in Fig.6, is executed. We set N equal to 5 and find the best value of E or so called $E_{min}$ corresponding to N=5. Suppose this value is $E_{min}$=1.0. Since this value is 1.0 and greater than an arbitrary set number of 0.1, we subtract 0.1 from 1.0 which results in 0.9. This is now considered as the designated error tolerance value for N=6 in the next round of iteration.
3) We continue based on the flowchart in Fig.6 with N=6 and E=0.9 as designated error tolerance value.
    a. If it does not reach the designated error tolerance of E=0.9 before 200 iterations for symmetrical cases and 300 iterations for nonsymmetrical cases (as shown in Fig.5), we move on to N=7 and E=0.9
    b. If it reaches the designated error tolerance of E=0.9 before 200 iterations for symmetrical cases and 300 iterations for nonsymmetrical cases (as shown in Fig.5); then, the program in Fig.6, is executed. We set N equal to 6 and find the best value of E or so called $E_{min}$ corresponding to N=6. Suppose this value is found to be $E_{min}$=0.01. Since this value (0.01) is less than 0.1, we divide 0.01 by 10 which results in 0.001. This is now considered as the designated error tolerance value for N=7 in the next round of iteration.
4) The mapping process gets repeated for all N ranging from 5 to 100 according to the flowchart in Fig.7. The last value of E and its corresponding N for which flowchart in

Fig.6 is executed, are the sought optimum values of E and N for the considered geometrical section.

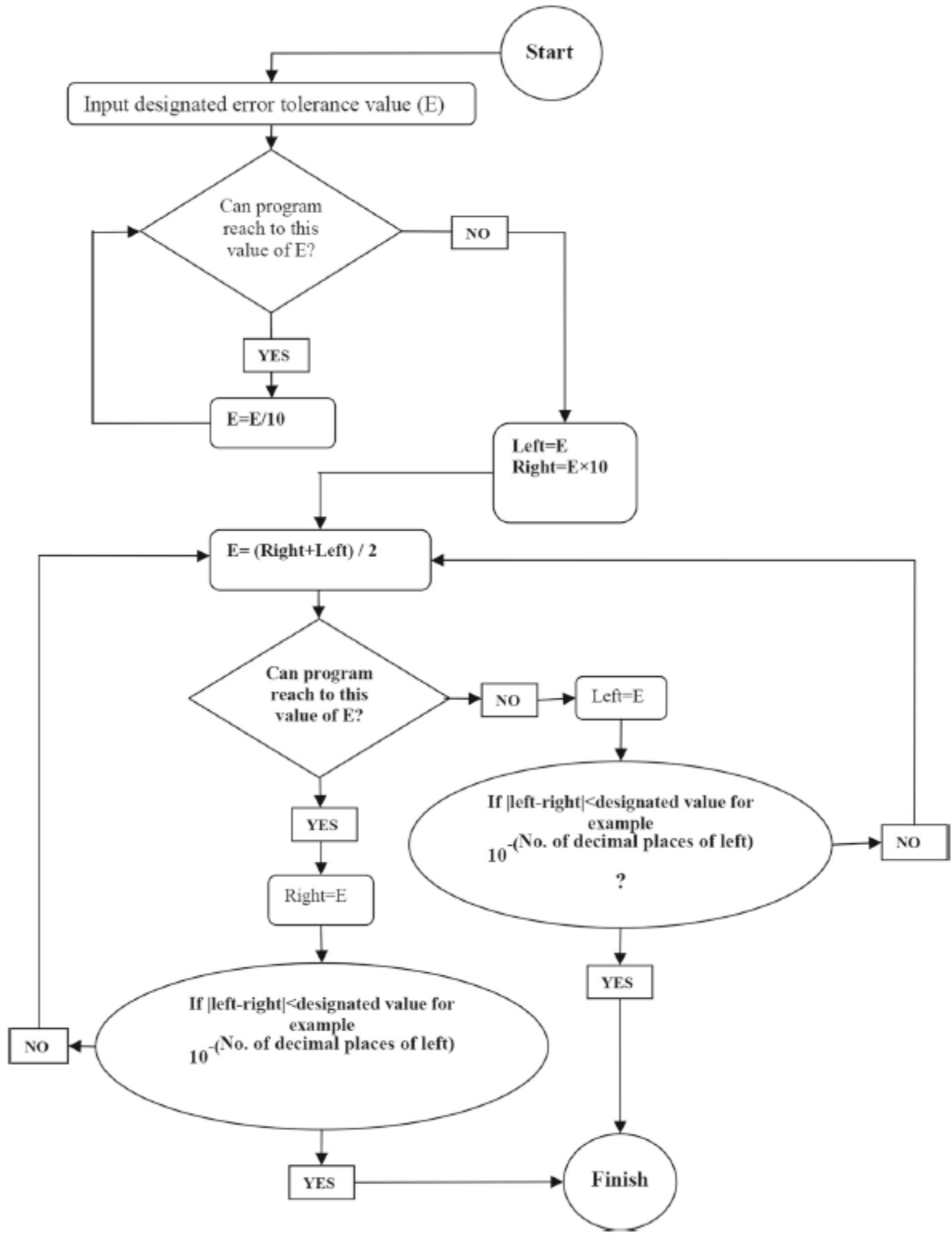

Fig.6. The flowchart of finding the best values of E for each N.

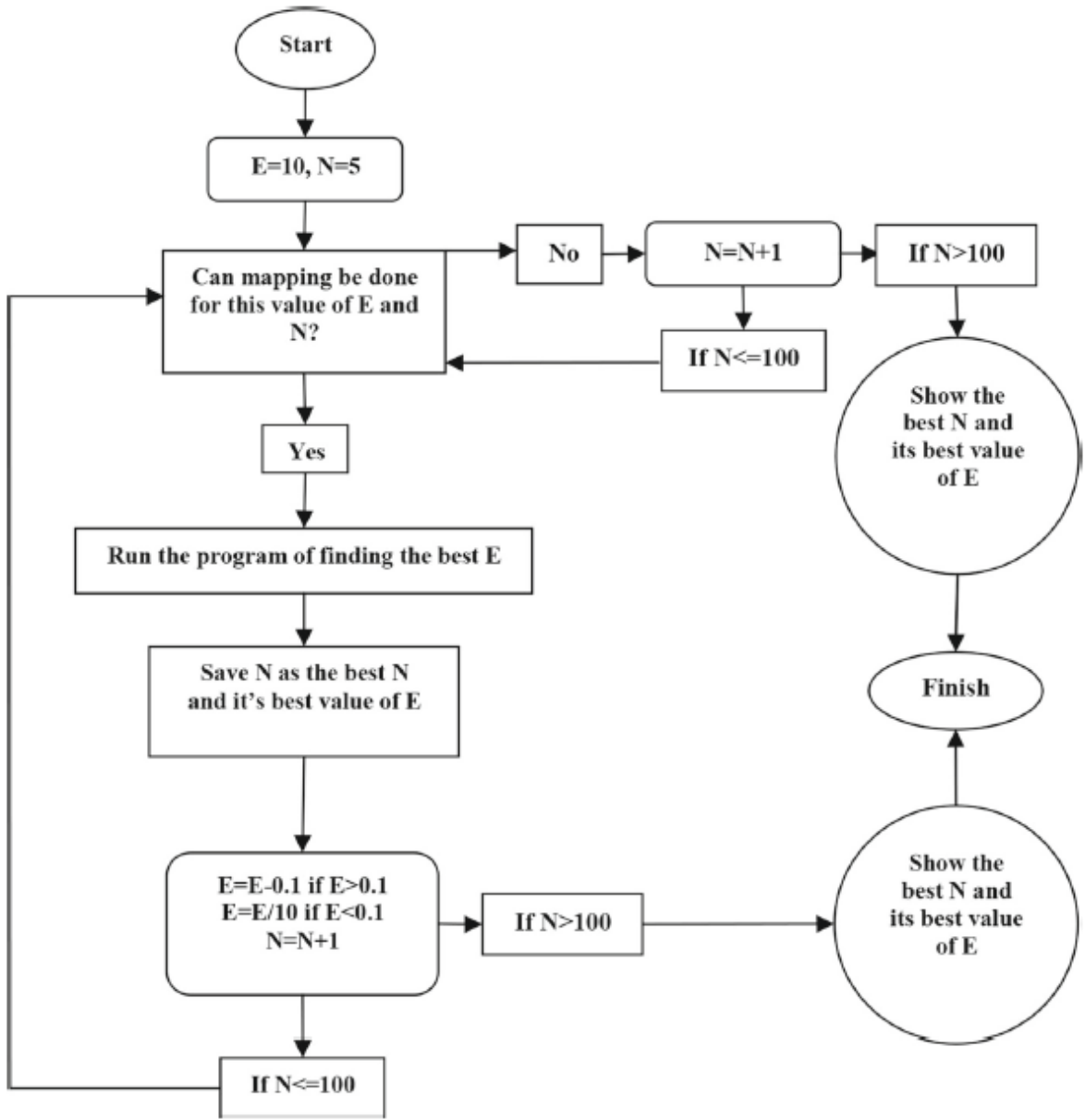

Fig.7. The flowchart of finding the best values of E and N in the mapping process.

## 4. Examples and discussions

The results of mapping different sections are presented by showing the best values of E that has been obtained from Eq.18 for different values of number of mapping coefficients (N).

These values of E have been based on N and are the smallest possible values of E for which a convergent solution is not attainable. In other words, the most accurate mapping profile is achieved by this value of E.

Sections that have been presented as examples in this section include rectangular section, large bulbous section, fine section and chine section. In each of these examples, the results of mapping the symmetrical sections, the results of mapping non-symmetrical sections, and their comparison against the results of Westlake and Wilson [23] have been presented, respectively. However, it must be noted that comparison of the current results in the case of chine sections with that of Westlake and Wilson [23] has not been presented because their study did not include this particular section.

For comparison of the obtained results for the symmetrical sections against that of the Lewis mapping, they are plotted next to the points of the real sections in Figs.8 (a), 10 (a), 14 (a), and 15 (a). As evidenced in these figures, the produced mapped points by the current method perfectly match the points of the real section, while the points produced by the Lewis mapping have considerable difference with the real section points.

For a better view of the relations between the mapped section points and their corresponding angles in the unit circle plane, the particular case of symmetric rectangular section and its transformation in the unit circle plane are demonstrated in Fig.9. In this figure, some points of the real section and their corresponding angles in the unit circle plane have been marked with the same numbering (or identifier).

As mentioned earlier, to demonstrate the capabilities of the current mapping method in the cases of rectangular, large bulbous and fine section for both symmetrical and non-symmetrical form, comparison of current mapping results and Westlake and Wilson Mapping [23] results is presented. The results indicate that in symmetrical cases, current mapping method is more accurate and robust, while in non-symmetrical cases, the results are similar to that by Westlake and Wilson mapping method [23].

### 4-1. Rectangular Section

Cross section of most of the barges, the floating breakwaters and the floating pontoons is rectangular. Because of hard chine in the corners, the mapping of these sections has more difficulty than the common sections. In other words, the mapped section and the real section will not agree with each other. However, mapping of these sections by the present method has produced very good results and shows high compatibility with the real sections. In this section, two examples are presented; one for the symmetrical and the other for the non-symmetrical rectangular section. The non-symmetrical section is a rectangular section at 15 degrees of heel.

To better understand the relation between the section points and their corresponding angles in unit circular plane, the section points and their corresponding angles are displayed in Fig.8 for the case of symmetric rectangular section.

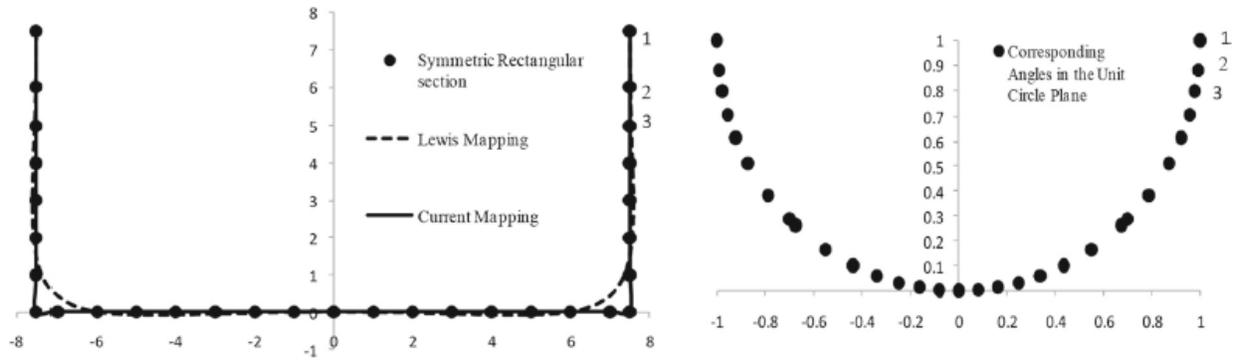

Fig.8. Illustration of a Symmetric rectangular section (left) and its corresponding angles in the transformed unit circle plane (right).

The results of mapping for the symmetrical and non-symmetrical rectangular sections are shown in Figs.9 and 11. The result of Lewis mapping technique is also demonstrated in Fig.9 for the case of symmetrical section. The Minimum values of E which are obtained for different values of N in mapping of symmetric rectangular section and non-symmetrical section have been shown in Figs.10 and 12

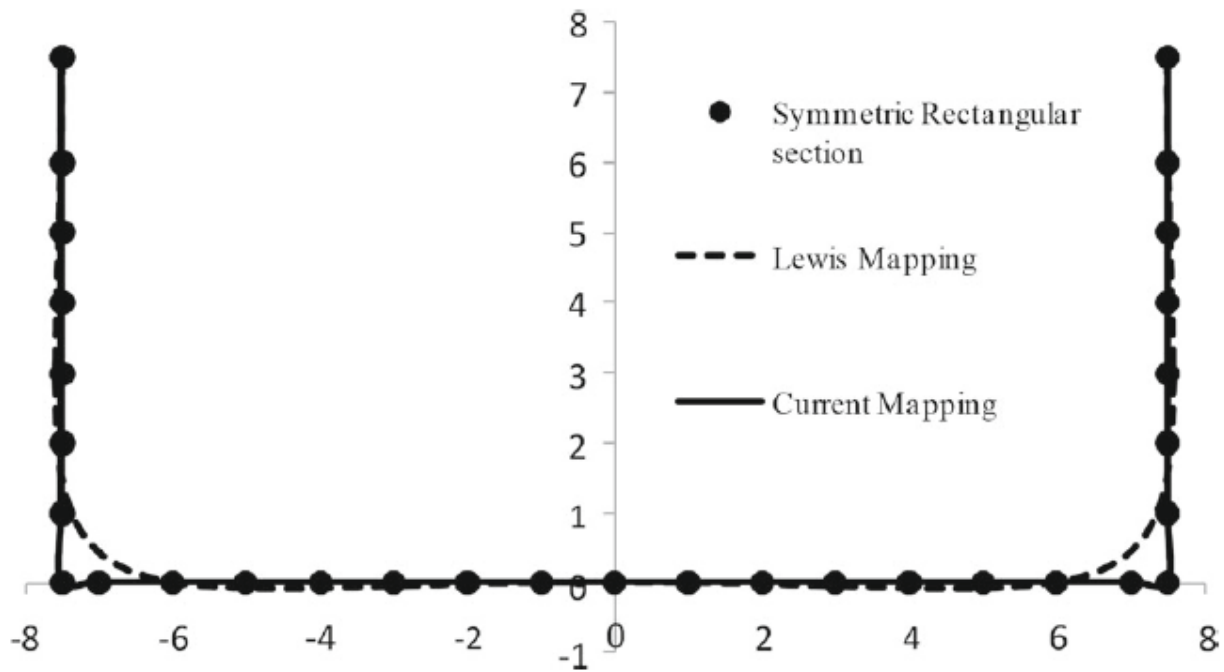

Fig.9. Results of conformal mappings for symmetric rectangular section.

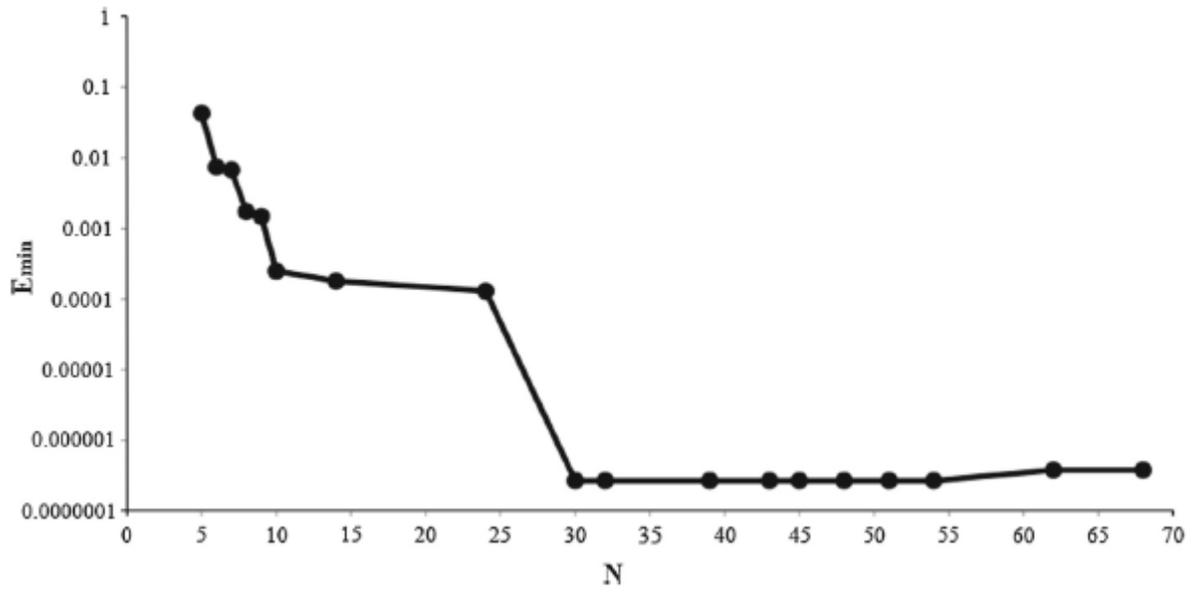

Fig.10 . Minimum value of E which is obtained for different values of N in mapping of symmetric rectangular section.

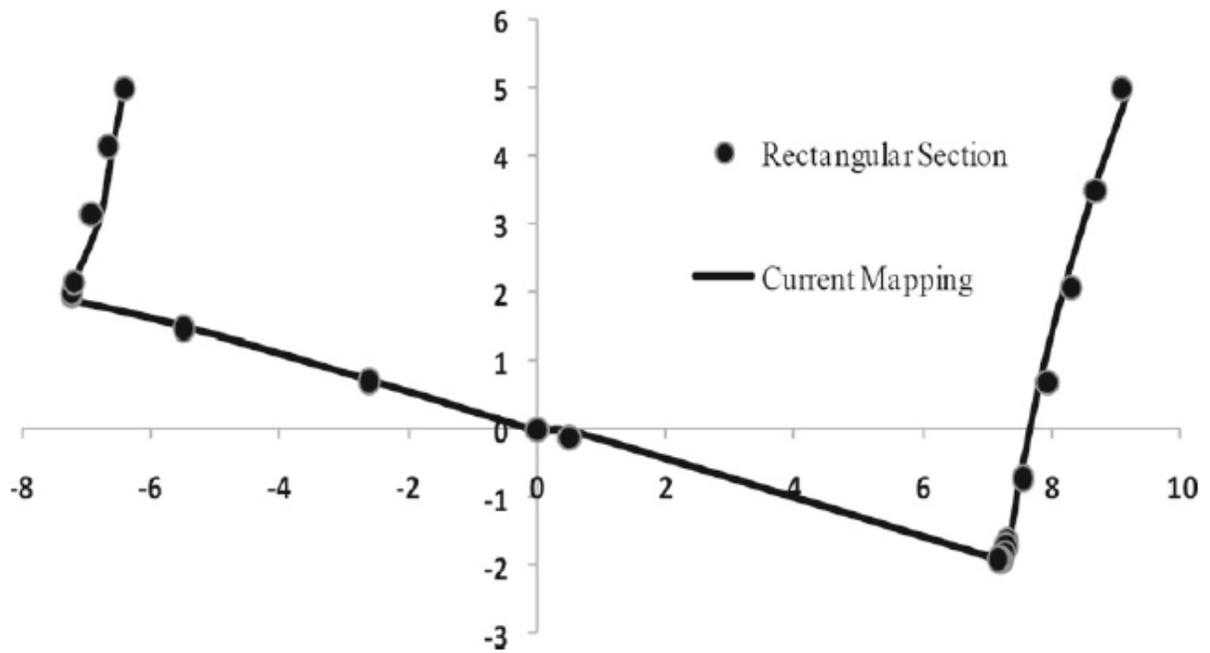

Fig.11 Results of conformal mappings for non-symmetric rectangular section.

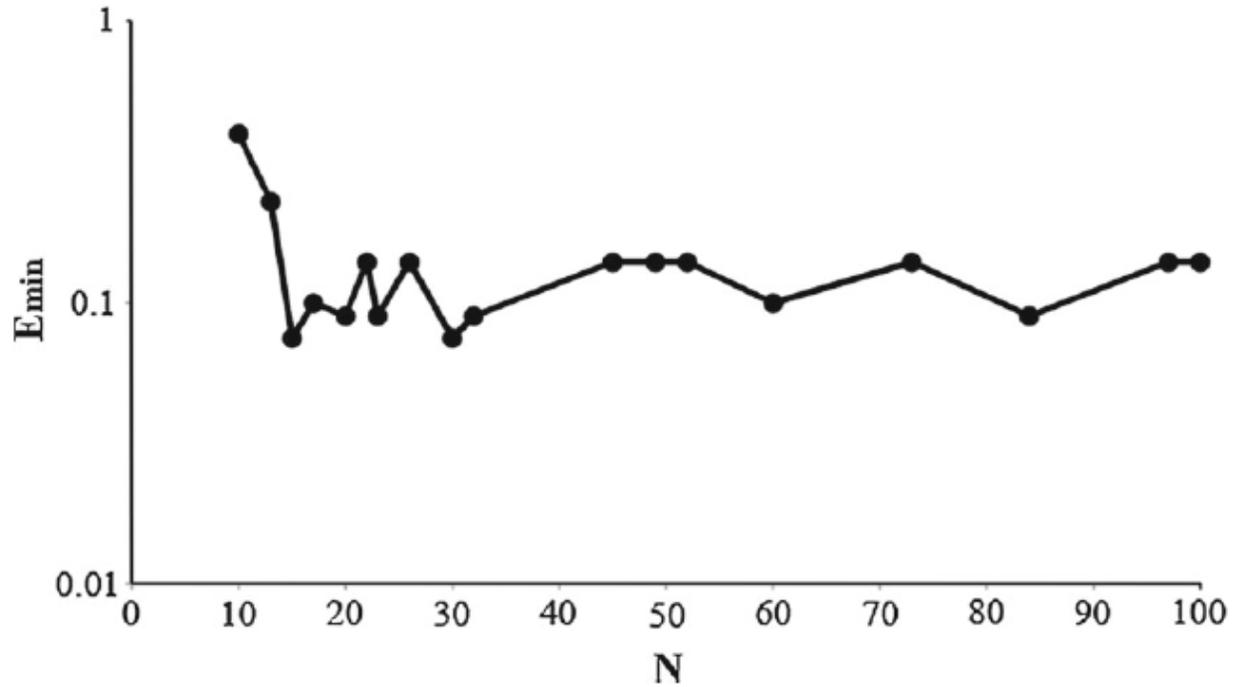

Fig.12. Minimum values of E which are obtained for different values of N in mapping of non-symmetric rectangular section.

For the purpose of comparing the current results against those by Westlake and Wilson [23], the current scheme has been applied to the same examples which are solved in Westlake and Wilson's article [23]. Accordingly, accuracy of the proposed method has been examined for two symmetrical and non-symmetrical sections. Coordinates of the points related to the examined sections along with the points resulting from the Westlake and Wilson's mapping [23] and their final mapped points are respectively plotted in Figs.13 and 14.

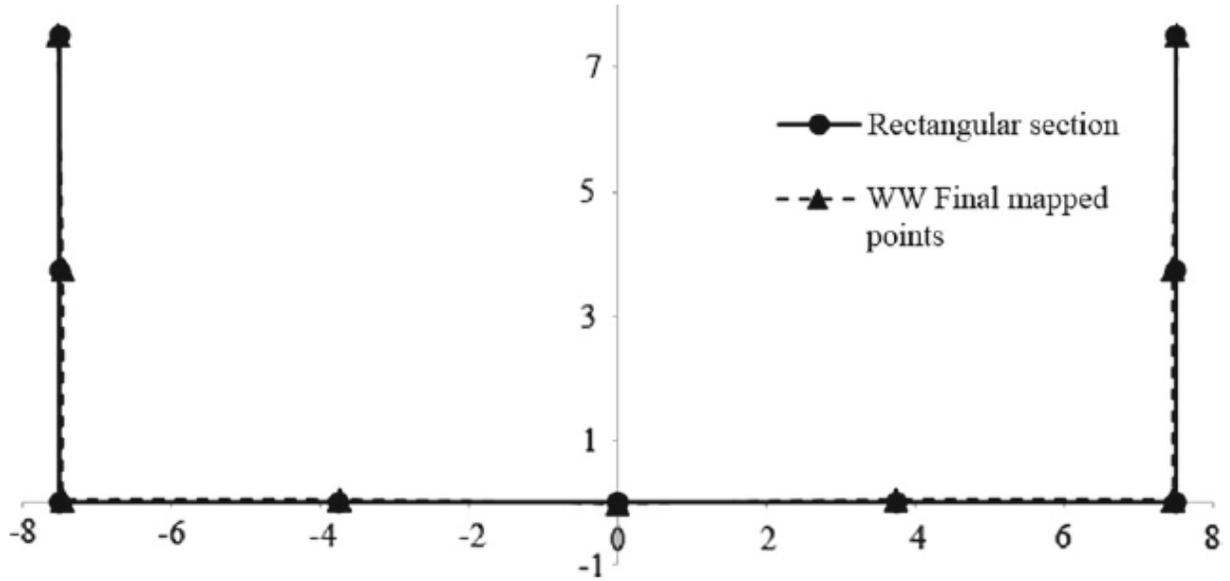

Fig.13. WW [23] offsets and final mapped point by using 12 transform parameters for symmetric rectangular section.

Therefore, error from the results of Westlake and Wilson's mapping [23] for N=12 is equal to

$$E_{WW} = \sum_{i=1}^{I}(x_i - x_{0i})^2 + (y_i - y_{0i})^2 = 0.006327$$

Results of the proposed technique for the current section are given in table 1, which for N=12 is equal to 3.0E-9. This indicates greater accuracy compared to Westlake and Wilson's [23] results.

Table 1. Minimum values of E that are obtained for different values of N by using the current mapping of symmetric rectangular section.

| N | 5 | 6 | 12 |
|---|---|---|---|
| E | 8.0 E-3 | 178.0 E-8 | **3.0E-9** |

It is quite apparent that, in order to reach the same level of accuracy as in the Westlake and Wilson's [23] mapping with N=12 parameters, it is sufficient to use N=6 parameters in the proposed method, which is indicative of robustness and fastness of the current technique.

Next, comparison is made for the case of a non-symmetric rectangular section with heel angle of 15 degrees which is examined in the work by Westlake and Wilson [23]. Geometry of the section and final mapped points is plotted in Fig.14.

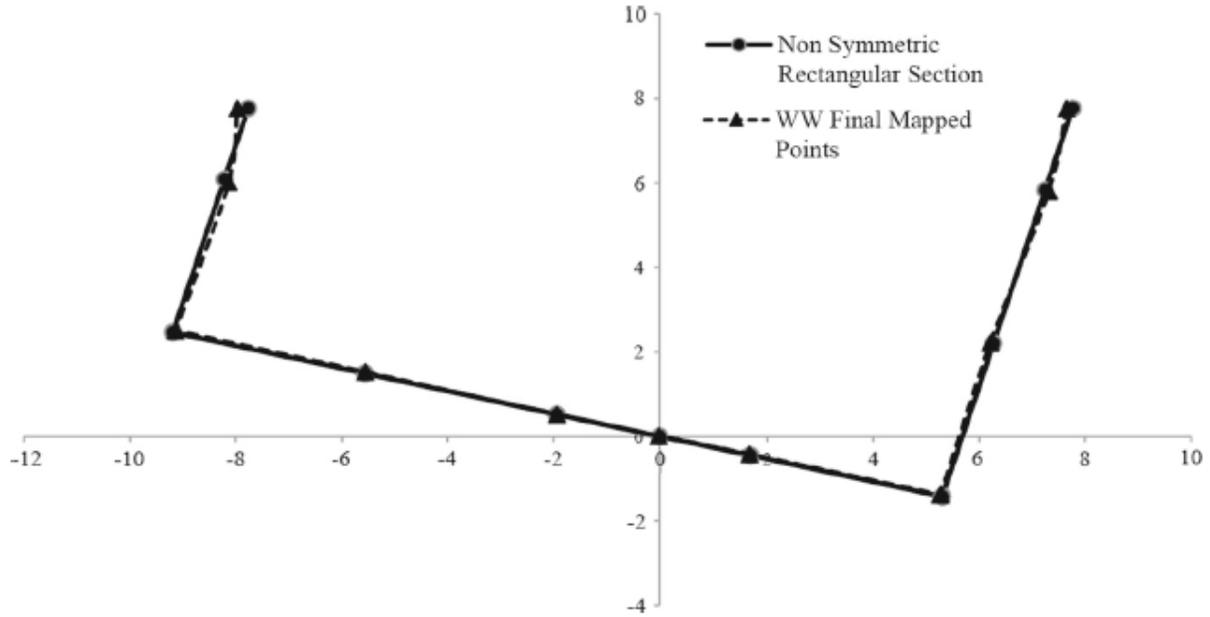

Fig.14. WW [23] offsets and final mapped point by using 12 transform parameters for non- symmetric rectangular section at 15 degrees of heel.

Here, error from the results of Westlake and Wilson's mapping [23] for N=12 is equal to

$$E_{WW} = \sum_{i=1}^{I}(x_i - x_{0i})^2 + (y_i - y_{0i})^2 = 0.084522$$

Results of the proposed technique for the current section are given in table 2.

Table 2. Minimum value of E which is obtained for different values of N by using current mapping for the non-symmetric rectangular section at 15 degree of heel.

| N | 12 | 15 | 17 | 30 |
|---|---|---|---|---|
| E | 0.115 | 0.079 | 0.06 | **0.004** |

It is clear that, error for N=12 in the current case is almost the same as the error from Westlake and Wilson's [23] resutls and that, in order to increase the accuracy in comparison with WW's results, N has to be slightly increased further, i.e. N=15.

### 4-2. Bulbous sections

Almost all of the more novel commercial ships have bulbous sections at their bows. Sections with large bulb can not be mapped easily. The two examples that are discussed here for the symmetrical as well as non-symmetrical bulbous sections indicate the great ability of the present mapping method for mapping such large bulbous sections. The non-symmetrical section is a bulbous section at 15 degrees of heel.

The section points and results of current mapping technique for symmetric and non-symmetric mapping are illustrated in Figs.15 and 17, while the errors are plotted in Figs.16 and 18. The result of Lewis mapping technique is also included in Fig.15 for the case of symmetrical section.

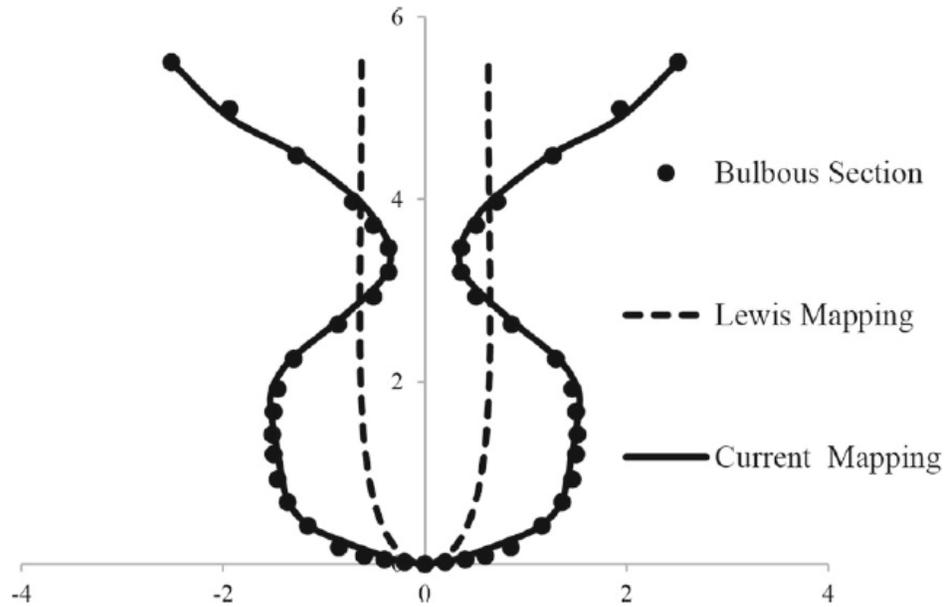

Fig.15. Results of conformal mappings for symmetric bulbous.

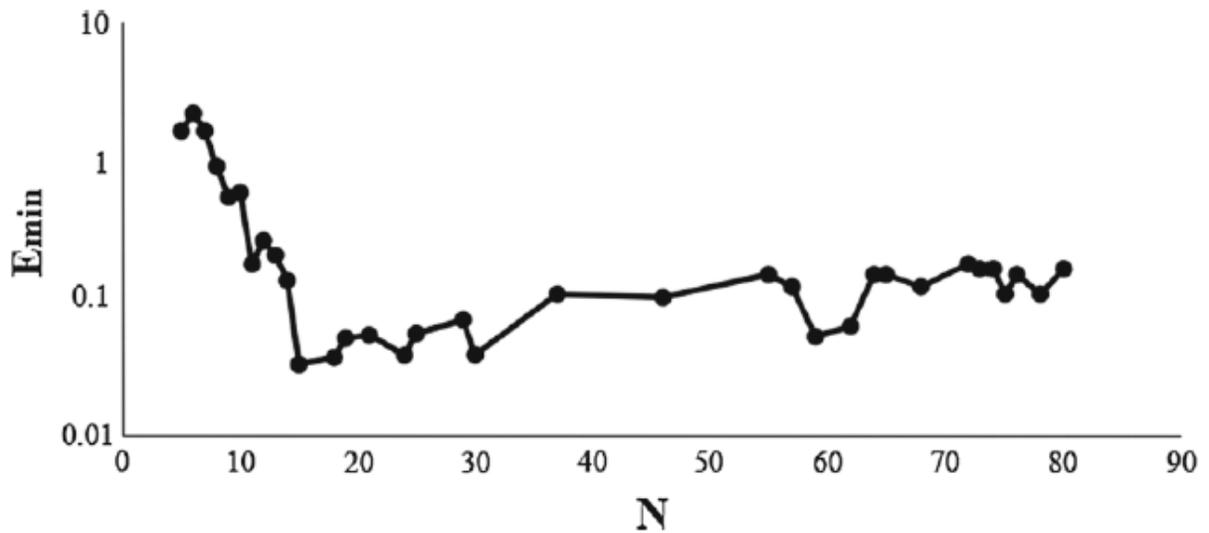

Fig.16. Minimum values of E which can be obtained for different values of N in mapping of symmetric large bulb section.

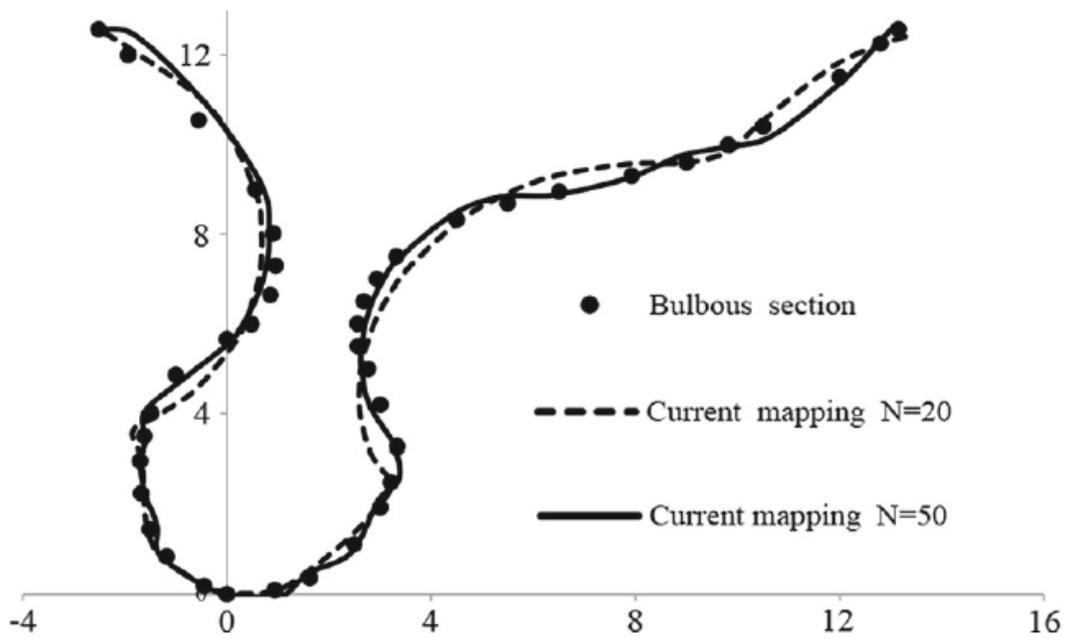

Fig.17. Results of conformal mappings for non-symmetric bulbous section.

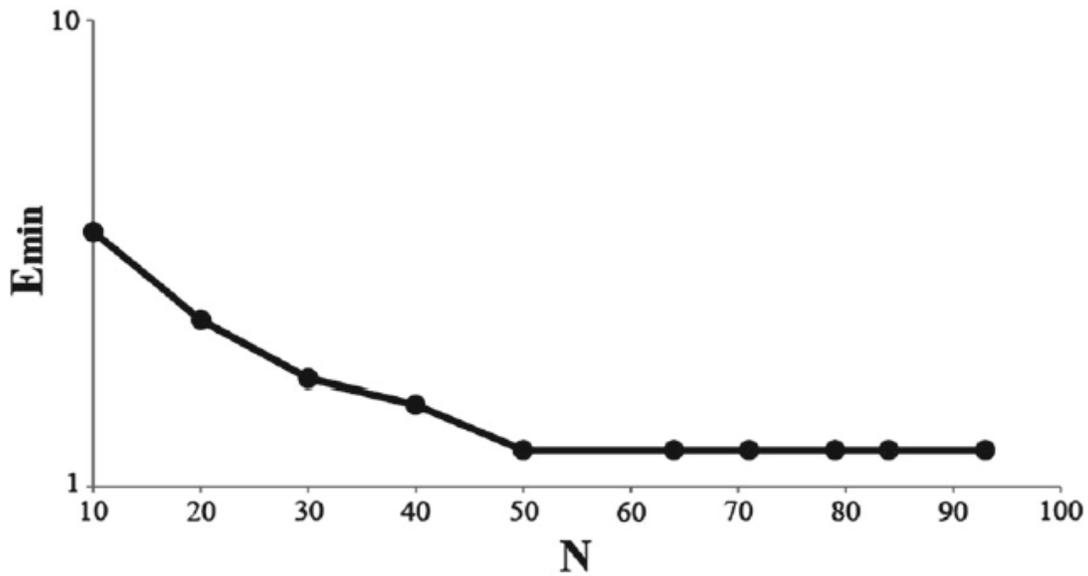

Fig.18. Minimum values of E which are obtained for different values of N in mapping of non-symmetric Bulbous sections at 15° of heel.

Results of the current mapping technique for the symmetric Bubous section (in the cases of N=19 and N=31) are illustrated and compared with the results of the Westlake and Wilson [23] (in the case of N=48) in Fig.19. As evidenced in these figures, the proposed technique offers more accurate results (compared to the real sections) with a fewer mapping parameters than the Westlake Wilson's [23] results.

The mapping of non-symmetric Bulbous section with 15 degrees heel, which was also done by Westlake and Wilson [23], has been done by the current mapping technique, too. Results of mapping of this type of section for N=16, 20, 28, 36, and 48 parameters have been graphically presented in Fig.20. In these plots, the currents results as well as those by Westlake and Wilson [23] are illustrated and compared. As evidenced in these plots, there is not much difference between the results of the two techniques and both methods posses similar error at various N. Therefore, one can easily conclude that the current mapping technique, eventhough offers compatible results with those by Westlake and Wilson [23] in the case of non-symmetric mapping bulbous section, but in the case of symmetric mapping, offers more accurate results in a shorter time period.

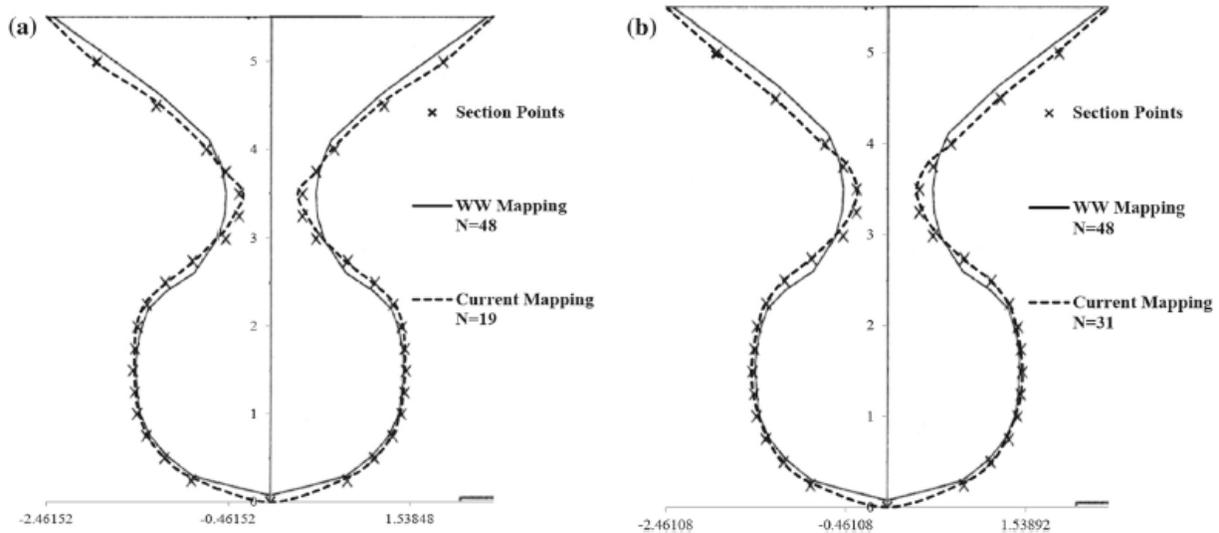

Fig.19. Results of conformal mappings for symmetric bulbous section for   a) N=19, compared with the results of WW [23] mapping for N=48   b) N=31, compared with the results of WW [23] mapping for N=48.

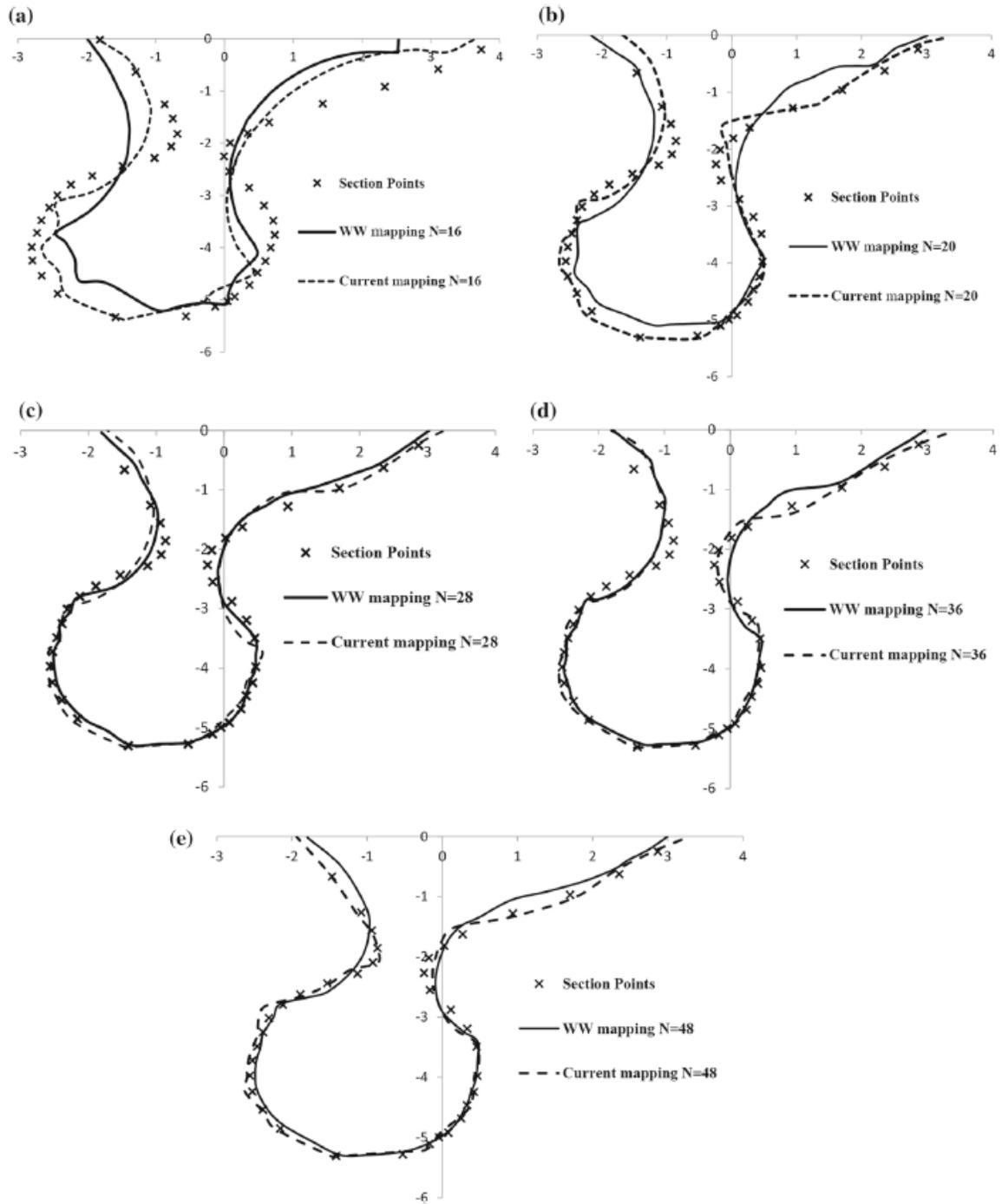

Fig.20. Comparison of the current results and that of WW [23] mapping for non-symmetric bulbous section for a) N=16, b) N=20, c) N=28, d) N=36, and e) N=48.

## 4-3. Fine section

Mapping of the fine sections is also a cumbersome task and most of the mapping methods do not have the ability to produce the mapped sections with good accuracy. The aft and fore sections of the ship, frigate and destroyer ship sections are considered as fine sections. The following examples for symmetrical and non-symmetrical sections illustrate the proper accuracy of results of the current mapping method. The non-symmetrical section is a fine section at 20 degree of heel.

Results of mapping a symmetric and non-symmetric fine sections are demonstrated in Figs.21 and 23. In this figure, the results of Lewis mapping for the case of symmetric mapping is also presented. The resulting errors for these mappings have been displayed in Figs.22 and 24.

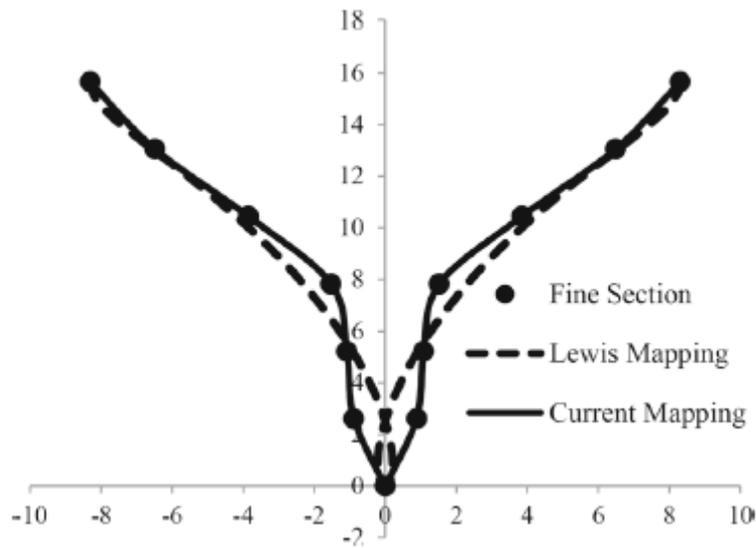

Fig.21. Results of conformal mappings for the symmetric fine section.

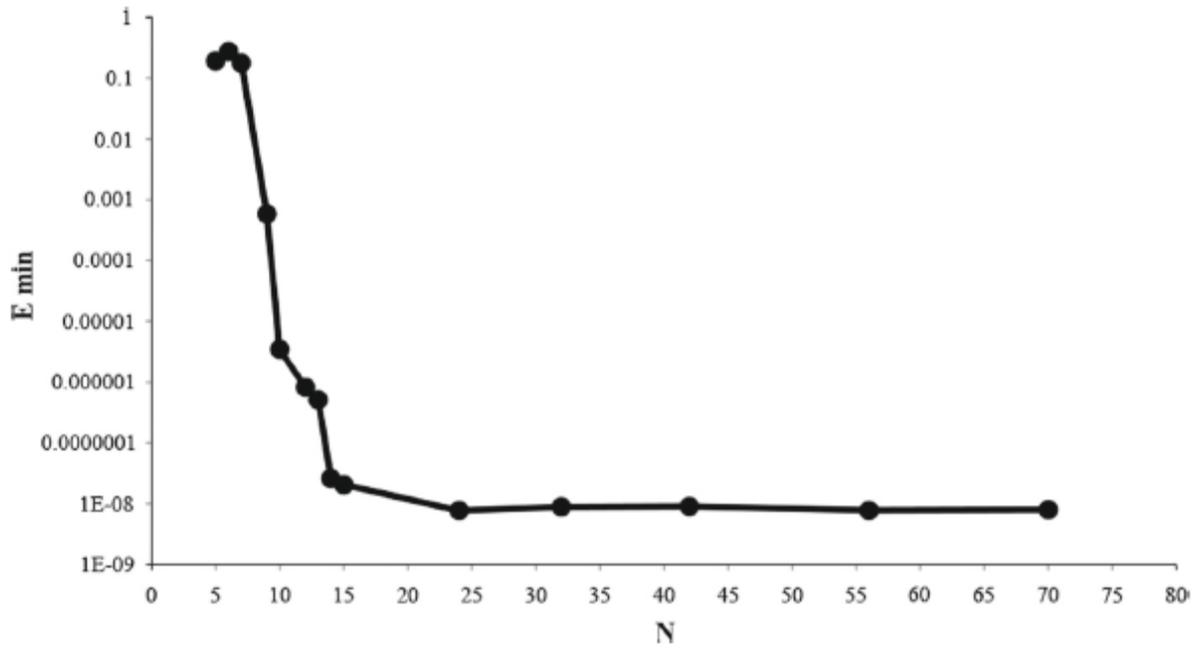

Fig.22. Minimum values of E which are obtained for different values of N in mapping of symmetric fine section.

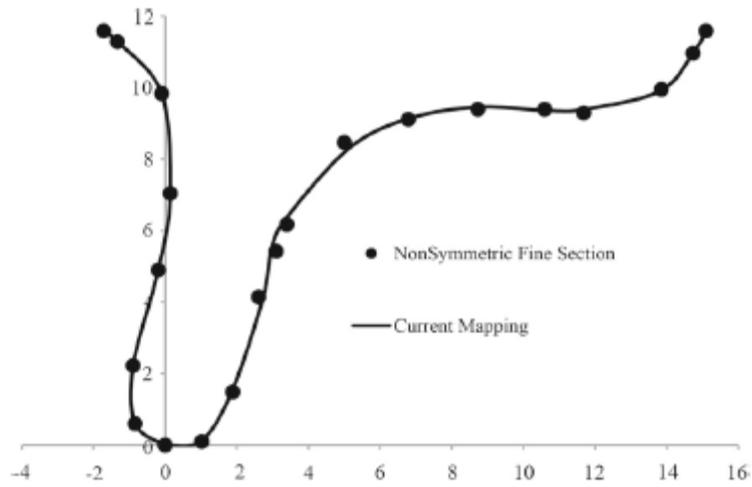

Fig.23. Results of conformal mappings for non-symmetric fine section.

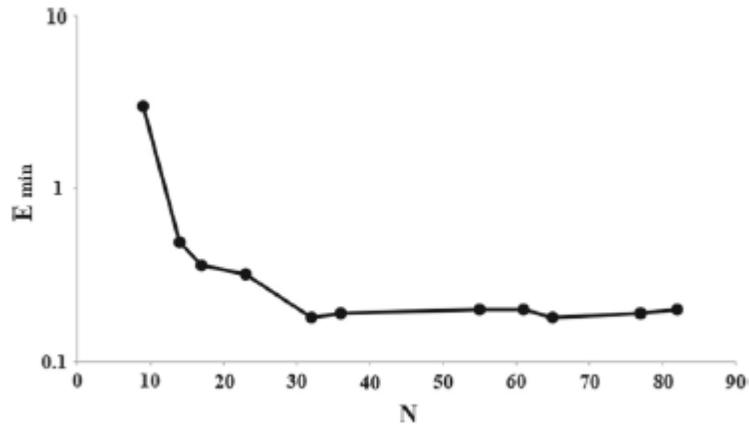

Fig.24. Minimum values of E which can be obtained for different values of N in mapping of non-symmetric fine section.

Coordinates of the symmetric and non-symmetric fine sections analyzed earlier by WW [23] are plotted in Figs.25 and 27.

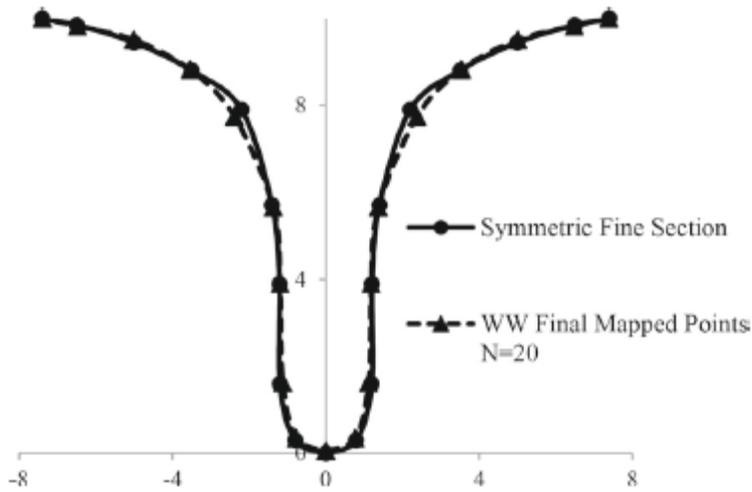

Fig.25. WW [23] offsets and final mapped point by using 20 transform parameters for symmetric fine section.

Error resulting from WW [23] method for N=20 is equal to

$$E_{WW} = \sum_{i=1}^{I}(x_i - x_{0i})^2 + (y_i - y_{0i})^2 = 0.081552$$

The results of the current mapping technique for symmetric fine section have been shown in Fig.26.

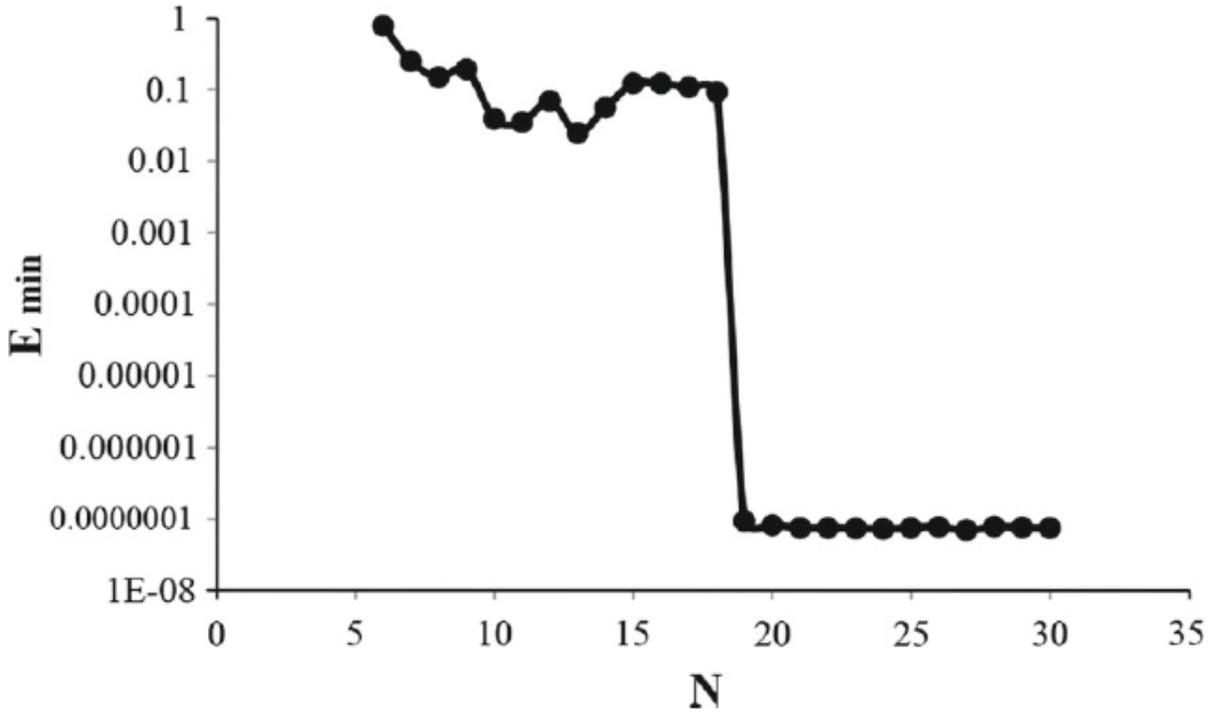

Fig.26. Minimum values of E which can be obtained for different values of N in mapping of symmetric fine section.

In the case of non-symmetric fine section, coordinates of the section and the results of WW [23] technique are given in Fig.27.

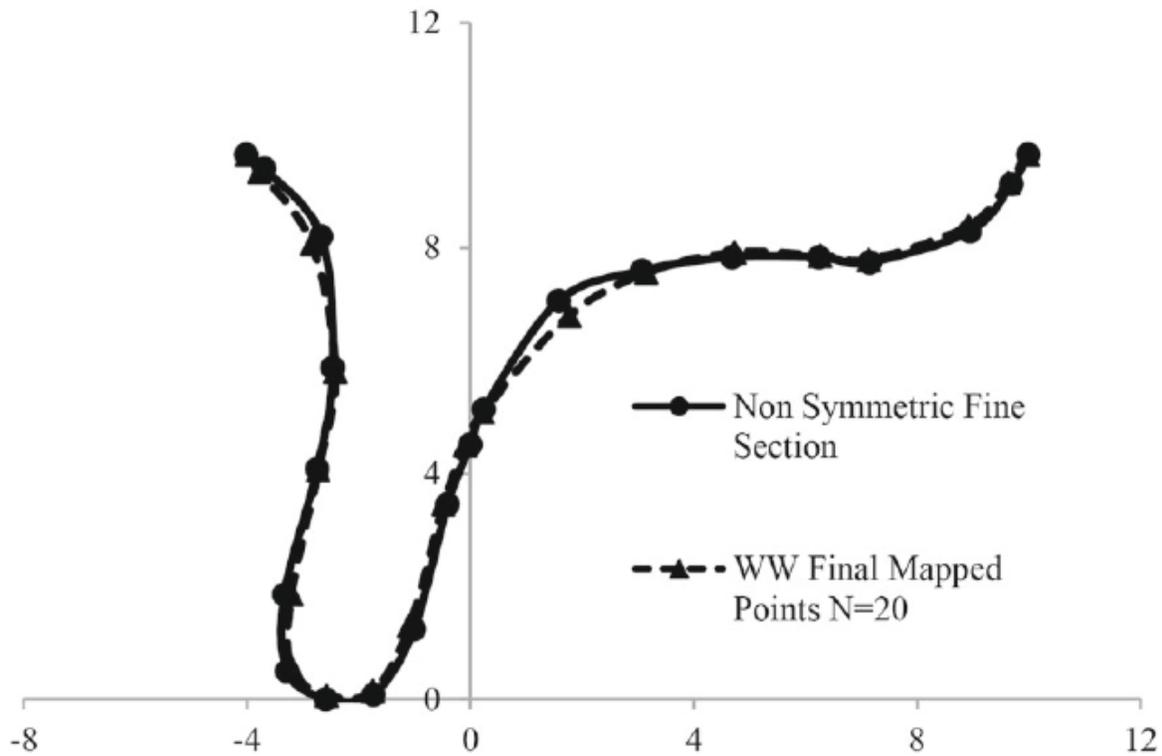

Fig.27. WW [23] offsets and final mapped point by using 20 transform parameters for non symmetric fine section.

Error resulting from the WW [23] mapping for N=20 is equal to

$$E_{WW} = \sum_{i=1}^{I}(x_i - x_{0i})^2 + (y_i - y_{0i})^2 = 0.267901$$

Results of the current method for the nonsymmetric fine section are tabulated in table 3.

Table 3. Minimum value of E which is obtained for different values of N in mapping of non-symmetric fine section.

| N | 16 | 18 | 19 | 20 | 26 | 30 | 35 |
|---|----|----|----|----|----|----|----|
| E | 0.19 | 0.2 | 0.256 | 0.21 | 0.18 | 0.14 | **0.13** |

It is clear that there is not much difference betwee the current results and that those by WW [23] technique for non-symmetric mapping of fine section.

### 4-4. Chine sections

Planning boats, catamarans and naval ships have chine to decrease the resistance and power consumption. Dynamic analysis of these vessels plays the main role in the design states. The following examples demonstrate the great mapping ability of the current method for the fine sections. The non-symmetrical example is a chine section at $20^0$ of heel angle.

Results from the symmetric as well as non-symmetric chine sections are illustrated in Figs.28 and 30. Also, the result of Lewis mapping for the case of symmetric mapping is presented in Fig.28.

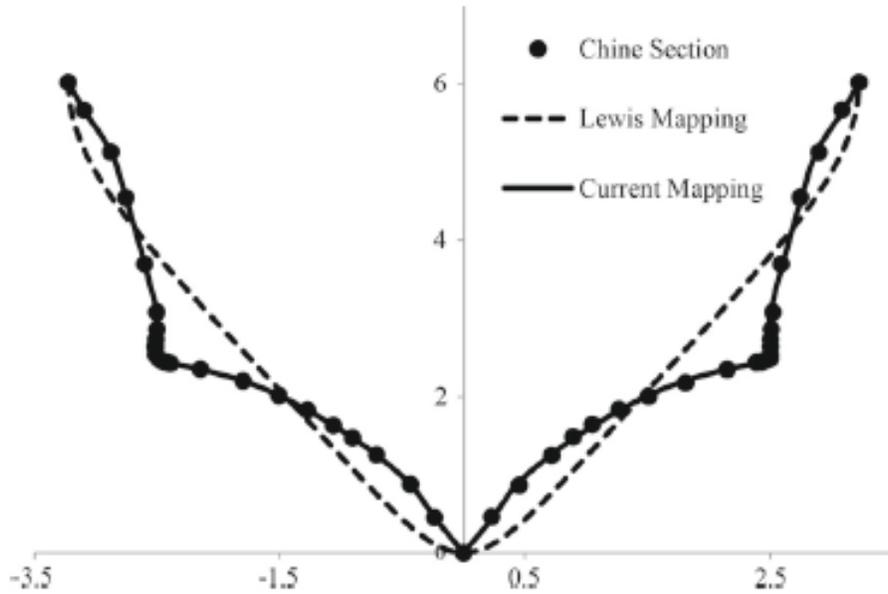

Fig.28. Results of conformal mappings for symmetric chine.

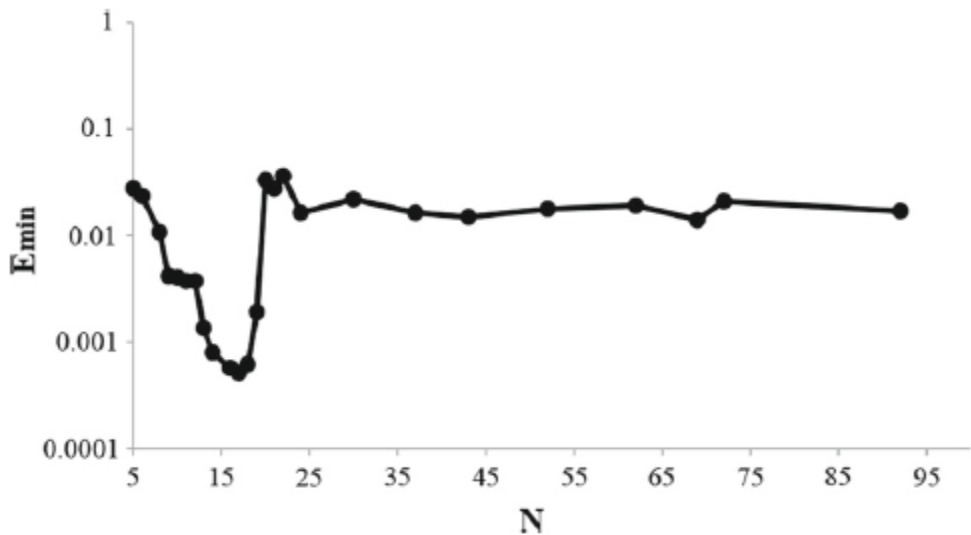

Fig.29. Minimum values of E which are obtained for different values of N in mapping of symmetric chine section.

Coordinates of the non-symmetric chine section and result of current mapping are plotted in Fig.30, while the results of mapping are presented in Fig.31.

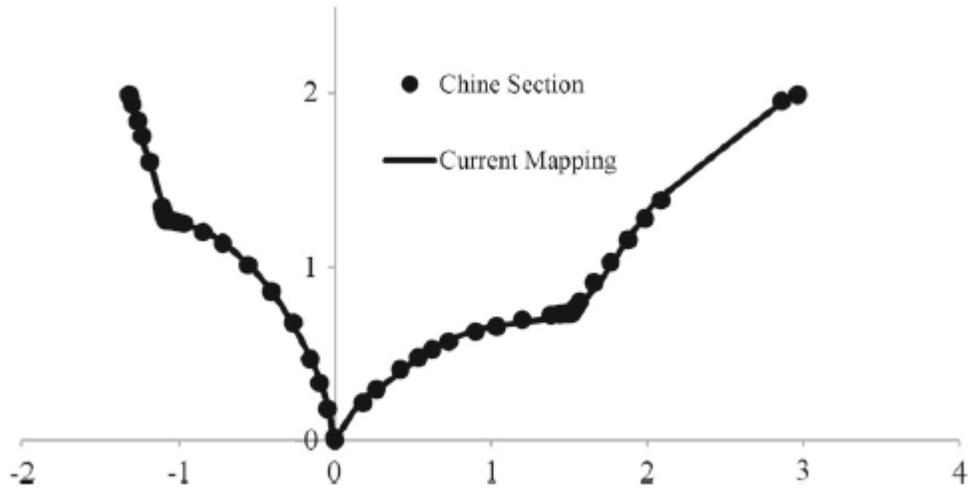

Fig.30. Results of conformal mappings for non-symmetric chine section.

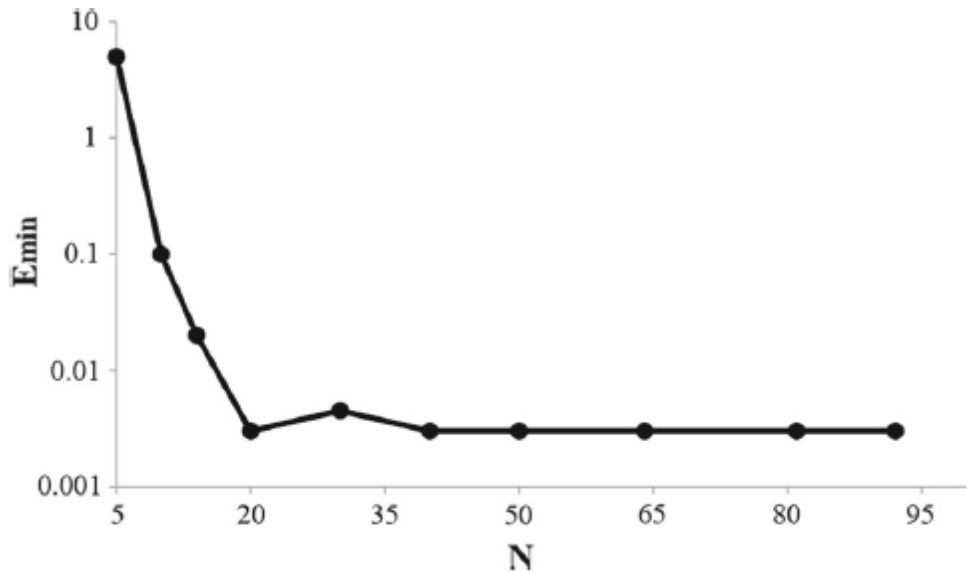

Fig.31. Minimum value of E which are obtained for different values of N in mapping of non- symmetric chine section.

Results of the current mapping technique for symmetric and non-symmetric chine section deomonstrated in Figs.28 and 30 indicate great match with the real section and suitable accuracy of the scheme.

### 4-5. Computational Times and Accuracy

Results of current mapping method were listed earlier for each of the geometrical sections. In this section, in addition to presenting the optimum values of error (E) and number of mapping coefficients (N) for each of these sections, computational times for finding these values. Also, effort has been made to show the accuracy associated with the proposed technique for determining these quantities. This is accomplished by using the Nash-Sutcliffe model efficiency coefficient (E).

### 4-5-1. Nash-Sutcliffe model efficiency coefficient

The Nash-Sutcliffe model efficiency or accuracy coefficient (E) is commonly used to assess the predictive power of hydrological discharge models [25]. However, it can also be used to quantitatively describe the accuracy of outputs for other physical models. Here, It is utilized for assessing the accuracy of the applied mapping and is defined as

$$E_x = 1 - \frac{\sum_{i=1}^{I}(X_{re,i} - X_{ma,i})^2}{\sum_{i=1}^{I}(X_{re,i} - \overline{X}_{re})^2}$$

$$E_y = 1 - \frac{\sum_{i=1}^{I}(Y_{re,i} - Y_{ma,i})^2}{\sum_{i=1}^{I}(Y_{re,i} - \overline{Y}_{re})^2}$$

where $E_x$ is the Nash-Sutcliffe model efficiency coefficient in the x-direction (for assessing the x-coordinate of the mapped points) and $E_y$ is the Nash-Sutcliffe model efficiency coefficient in the y-direction (for assessing the x-coordinate of the mapped points). Also, $X_{re,i}$ is the x-coordinate of the real section points, $X_{ma,i}$ is the x-coordinate of the mapped points, $\overline{X}_{re}$ is the mean of x-coordinates of the real section points, $Y_{re,i}$ is the y-coordinate of the real section points, $Y_{ma,i}$ is the y-coordinate of the mapped points, $\overline{Y}_{re}$ is the mean of y-coordinates of the real section points. Essentially, the closer the efficiency coefficient is to 1, the more accurate the mapping is. Table 4 displays the best values of error (E), number of mapping coefficients (N), the required time for finding the optimum values of E and N and the corresponding values of Nash-Sutcliffe model efficiency coefficients. As evidenced in this table, great accuracy has been achieved in the applied mapping technique for all the symmetric and non-symmetric sections, since the efficiency of mapping for both x and y coordinates of points ranges from 0.997217924 to 1.

Table 4. Optimum values of E and N, the required computational time, and the mapping efficiency.

| Sections | The best value of E | The best value of N | Log $E_{min}$ | Time to finding the best N and E(sec) | $E_x$ | $E_y$ |
|---|---|---|---|---|---|---|
| Symmetric rectangular | 2.6875E-07 | 30 | -6.57065 | 10.6324005 | 1 | 1 |
| Non-symmetric rectangular | 0.075 | 15 | -1.12494 | 12.9324015 | 0.999363752 | **0.999903685** |
| symmetric bulbous | 0.033203 | 15 | -1.47882 | 11.9324015 | 0.999843 | **0.997217924** |
| Non-symmetric bulbous | 1.2 | 50 | 0.079181 | 48.438707 | 0.998803 | **0.997949** |
| Symmetric fine | 7.68E-09 | 24 | -8.11464 | 11.96342 | 1 | **1** |
| Non- symmetric fine | 0.18 | 32 | -0.74473 | 13.96342 | 0.999825 | **0.99983** |
| Symmetric chine | 0.033203 | 20 | -1.47882 | 14.315382 | 0.99953946 | **0.999178** |
| Non- symmetric chine | 0.003 | 20 | -2.52288 | 12.245617 | 0.999906 | **0.999985103** |

## 5. Conclusions

A new generation of multi-parameter conformal mapping method has been proposed. Using this method, it has been demonstrated that any arbitrary section such as large bulbous, fine, rectangular and chine sections and sections with abrupt curvatures can be easily mapped to a unit circle with high accuracy and minimum error. Contrary to the previous generation of multi-parameter mapping techniques, the suggested method is drastically different in detail and can map any arbitrary symmetrical and non–symmetrical sections onto a unit circle in a much quicker way. It uses least square scheme to compute the mapping coefficients, considers a different value for the first mapping coefficient, solves all the mapping parameters simultaneously, and by making a simple assumption, an equation is obtained which is solved for the corresponding angles in the mapped unit circle plane. Moreover, it makes use of the results of the symmetrical case for assigning the initial guess for the mapping coefficients in the case of non-symmetrical sections. In this case, the initial guess is assumed to be the average of the two mapping coefficients found based on the approach adopted for the symmetrical sections and for the right and left (half) distance from the centerline on the waterline. The newly introduced techniques make the multi-parameter mapping technique more robust, shorten the process of mapping and decrease the computing time.

Several examples of mapping the symmetrical as well as non-symmetrical sections have been demonstrated with high accuracy. Results of the proposed mapping technique, in the case of symmetrical sections, have been sketched and compared against the results of the Lewis mapping. In the case of symmetrical sections, comparison between the current mapping method and the method of Westlake and Wilson [23] illustrates better agreement of the current results with the points of the real sections. Another significant point to be mentioned is the fact that, in order to achieve a designated accuracy, the current method needs lesser N value that makes the proposed mapping technique more accurate and faster. On the other hand, in the case of non-symmetrical mapping, comparison of the results of the current mapping method and Westlake and Wilson scheme indicated similarity of the results or negligible differences. To better display the relations between the points on the real section and their corresponding angles in the mapped unit circle plane, the corresponding angles for the symmetric rectangular section have been plotted in a figure.

Based on the suggested multi-parameter mapping technique, a unique program has been produced that can be joined with any two dimensional potential problem solvers for the hydrodynamics analysis of flow around different marine structures.